\documentclass[amssymb,twoside,12pt,amscd,reqno]{amsart}

\usepackage{latexsym,amsfonts}
\usepackage{amsmath,amsthm}          
\usepackage{amssymb}                 
\usepackage{amscd}                   
\usepackage{euscript}                
\usepackage{epsfig}
\usepackage{psfig}

\theoremstyle{plain}
\numberwithin{equation}{section}
\newtheorem{thm}{Theorem}[section]

\newtheorem{cor}[thm]{Corollary}
\newtheorem{defn}[thm]{Definition}

\newtheorem{lemma}[thm]{Lemma}
\newtheorem{notn}[thm]{Notation}

\newtheorem{obs}[thm]{Observation}
\newtheorem{prop}[thm]{Proposition}
\newtheorem{prop_def}[thm]{Proposition/Definition}

\newtheorem{rmk}[thm]{Remark}
\newtheorem{rmk_def}[thm]{Remark/Definition}

\newcommand{\bi}{\begin{itemize}}
\newcommand{\ei}{\end{itemize}}
\newcommand{\bp}{\begin{proof}}
\newcommand{\ep}{\end{proof}}

\def\CC{\mathbb{C}}

\def\PP{\mathbb{P}}

\def\ZZ{\mathbb{Z}}
\def\al{\alpha}
\def\be{\beta}
\def\de{\delta}

\def\ka{\kappa}

\def\om{\omega}
\def\si{\sigma}
\def\th{\theta}
\def\ups{\upsilon}

\def\ze{\zeta}

\def\Ga{\Gamma}

\def\Th{\Theta}

\def\A{\mathcal{A}}
\def\AA{\tilde{\A}}

\def\D{\mathcal{D}}

\def\Lb{\mathcal{L}}
\def\E{\mathcal{E}}
\def\F{\mathcal{F}}

\def\G{\mathcal{G}}

\def\HH{\mathcal{H}}
\def\I{\mathcal{I}}
\def\J{\mathcal{J}}

\def\K{\mathcal{K}}
\def\M{\mathcal{M}}
\def\MM{\overline{M}}
\def\N{\mathcal{N}}
\def\pp{\overline{p}}

\def\S{\mathcal{S}}

\def\S{\mathcal{S}}
\def\T{\mathcal{T}}
\def\U{\mathcal{U}}
\def\V{\mathcal{V}}

\def\Bb{\mathfrak{B}}
\def\Bbe{\mathfrak{B}_E}
\def\BBb{\overline{\mathfrak{B}}}

\def\Mae{\overline{\M(a,e)}}
\def\ME{\overline{\M_E}}

\def\Ss{\mathfrak{S}}

\def\SSs{\overline{\mathfrak{S}}}

\def\ad{\text{ad}}
\def\AA{\text{A}}
\def\Aut{\text{Aut }}
\def\codim{\text{codim }}
\def\coker{\text{coker }}
\def\deg{\text{deg}}
\def\det{\text{det}}
\def\dim{\text{dim }}

\def\ev{\text{ev}}

\def\Ext{\text{Ext}}

\def\h{\text{h}}
\def\H{\text{H}}
\def\Hilb{\text{Hilb}}

\def\Hom{\text{Hom}}

\def\id{\text{id}}

\def\length{\text{length }}
\def\Pic{\text{Pic}}
\def\PGL{\text{PGL}}
\def\Proj{\text{Proj}}

\def\rk{\text{rk }}
\def\R{\text{R}}
\def\Sch{\text{Sch}}
\def\Sets{\text{Sets}}
\def\Spec{\text{Spec}}
\def\Sym{\text{Sym}}

\def\TTor{\mathcal{T}or}

\def\dra{\dashrightarrow}

\def\lra{\longrightarrow}

\def\ra{\rightarrow}
\def\bo{\boxtimes}
\def\da{\dagger}
\def\eset{\emptyset}

\def\setm{\setminus}

\def\ps{\vspace{4pt}}

\begin{document}

\title{Rational Families of Vector Bundles on Curves, II}
\author[Castravet]{Ana-Maria Castravet}
\address{Institute for Advanced Study \\ Princeton, NJ 08540}
\email{noni@ias.edu}
\thanks{This material is based upon work partially supported by the National
Science Foundation under agreement No. DMS-9729992.}

\begin{abstract}
This is a continuation of ``Rational Families of Vector Bundles on
Curves, I''.\\
Let $C$ be a smooth complex projective curve of genus $g\geq2$
and let $M$ be the moduli space of rank $2$, stable vector bundles
on $C$, with fixed determinant of degree $1$. We prove that there
is a one-to-one correspondence between rational curves on $M$ and
rank $2$ vector bundles on $\PP^1\times C$ with some fixed Chern
classes. For any $k\geq1$, we find all the irreducible components
of the space of rational curves  of degree $k$ on $M$ and their
maximal rationally connected quotients. 
\end{abstract}

\maketitle

\tableofcontents

\section{Introduction}

\subsection{Statement of results}

\

\ps

Let $C$ be a genus $g\geq2$ smooth projective curve over $\CC$.
Let $\xi$ be a degree $1$ line bundle on $C$ and let $M$ be the
moduli space of isomorphism classes of stable, rank $2$ vector
bundles $\E$ on C, with determinant
$\det(\E)=\wedge^2(\E)\cong\xi$. Then $M$ is
a smooth, projective, irreducible variety of  dimension $(3g-3)$.
It is known that $\Pic(M)\cong\ZZ$ \cite{DN} and let $\Th$ be the ample
generator. The canonical bundle $K_M$ is $\Th^{-2}$ \cite{R}. 

\ps

Our main goal in this paper is to give a complete
description of the irreducible components of the Hilbert scheme 
$\Hom_k(\PP^1,M)$ of rational curves on $M$ of 
degree $k\geq1$. By  a rational curve of degree $k$ on $M$ we will 
always mean a non-constant morphism $f:\PP^1\ra M$, such that the line 
bundle $f^*\Th$ on $\PP^1$ has degree $k$.
The expected dimension of $\Hom_k(\PP^1,M)$ is
\begin{equation*}
\h^0(f^*T_M)-\h^1(f^*T_M)=\deg(f^*T_M)+\dim M=2k+3g-3
\end{equation*}

For integers $a$ and $e$, we define the range (\ref{star}) to be
\begin{equation*}\tag{$\star$}\label{st}
\{(a,e)|\quad k\geq a>k/2,\quad\frac{k-a}{2a-k}\geq
e>0\}\cup\{(k,0)\}
\end{equation*}

\begin{thm}
For any pair of integers $(a,e)$ in the range (\ref{st}), there
are irreducible subvarieties $\M(a,e)\subset\Hom_k(\PP^1,M)$ such
that $\Mae$ is an irreducible component of $\Hom_k(\PP^1,M)$ if
and only if  $\dim\M(a,e)\geq(2k+3g-3)$.
\end{thm}

In \cite{C} we proved that there is a \emph{nice component}
$\M\subset\Hom_k(\PP^1,M)$: it has the expected dimension and the
general element is unobstructed. Moreover, the general element is
a very free curve if $k$ is sufficiently large. 

\begin{thm}
Let $\M\subset\Hom_k(\PP^1,M)$ be the nice component. If $k$ is
odd, then $\M$ is one of the varieties $\Mae$, namely the one for
$a=\frac{k+1}{2}$ and $e=\frac{k-1}{2}$. If $k$ is even, the nice
component $\M$ is not one of the varieties $\Mae$. For any $k$,
the nice component together with the irreducible components $\Mae$
are all the irreducible components of $\Hom_k(\PP^1,M)$.
\end{thm}

In \cite{C} we proved that if $g$ is even and  $k$ is divisible by
$(g-1)$, there is at least one extra component, which we called the 
\emph{almost nice component}: it has the expected dimension and the
general element is a free curve $f:\PP^1\ra M$ with the property that
$f^*T_M\cong O^g\oplus\A$, for some positive vector bundle $\A$.

\begin{thm}
If $g$ is even and  $k$ is divisible by $(g-1)$, let
$\M'\subset\Hom_k(\PP^1,M)$ be the almost nice component. Then
$\M'$ is the variety $\Mae$ for
$a=\frac{kg}{2(g-1)}, e=\frac{g}{2}-1$. The almost nice component $\M'$
together with the nice component $\M$ are the only irreducible components 
of the space $\Hom_k(\PP^1,M)$ that have the  general point unobstructed.
\end{thm}

Recall that a projective variety $X$ is 
\emph{rationally connected}, if 
there is a dense open $X^0\subset X$, such that for any two points $x_1$
and $x_2$ in $X^0$, there is a rational curve $\PP^1\ra X$ through $x_1$ 
and $x_2$. If $X$ is a smooth projective variety, its 
\emph{maximal rationally connected fibration (MRC)} 
is a rational map $\psi:X\dra Z$ such that: i) the general fiber is 
rationally connected, and ii) for very general $z\in Z$, a rational curve 
in $X$ intersecting the fiber $\psi^{-1}(z)$ is contained in $\psi^{-1}(z)$.
The MRC fibration of a smooth projective variety is known to exist \cite{K}.
The pair $(Z,\psi)$ is unique up to a birational transformation.

\ps

We are interested in the question of finding the MRC fibrations
of the spaces of rational curves on rationally connected varieties.
The following theorem answers this question for the case of 
moduli spaces $M$ of vector bundles.
\begin{thm}
Let $(a,e)$ be a pair of integers in the range (\ref{st}) and
let $\de=(k-a)-e(2a-k)$. We distinguish the following cases: \bi
\item[i) ] If $\de=0$, the MRC fibration of the variety $\M(a,e)$
is given by a surjective map $$\M(a,e)\ra J(C)$$ \item[ii) ]If
$\de>0$, the MRC fibration of the variety $\M(a,e)$ is given by a
map $$\M(a,e)\ra J(C)\times J(C)$$ The map is dominant if and only
if $\de\geq g$. \ei
\end{thm}

These theorems recover some of the results in \cite{C}. (Our methods
do not prove for example the result in \cite{C} about the general element
of the nice component being a very free rational curve if $k$ is sufficiently
large.) However, our proofs here rely on some of the constructions and results
in \cite{C}.

\subsection{Main Idea}

\

\ps

The key point in all of the above theorems is the following fact: to give
a morphism $\PP^1\ra M$ is the same as giving a rank $2$ vector bundle $\F$
on $\PP^1\times C$, satisfying the condition that the bundle
$\F_p=\F_{|\{p\}\times C}$  is stable for any $p\in\PP^1$.
Such a bundle $\F$ is unique up to tensoring
with a line bundle from $\PP^1$.

\ps

On $M\times C$, there is a rigidified {Poincar\'e} vector
bundle $\U$: if $c\in C$, then
$\det(\U_{|M\times\{c\}})\cong\Th$. It follows that associating to
a morphism $f:\PP^1\ra M$ the rank $2$ vector bundle $(f\times\id_C)^*\U$ on
$\PP^1\times C$, gives a one-to one correspondence between rational curves 
on $M$ and rank $2$ vector bundles on $\PP^1\times C$, satisfying the 
stability condition. 

\ps

As a consequence, we are lead to consider moduli 
spaces for rank $2$ vector bundles $\F$ on $\PP^1\times C$ with fixed 
Chern classes, constructed by Brosius in \cite{BR1} and \cite{BR2}, in the 
more general case of rank $2$ vector bundles on a ruled surface.
Precisely, we have to consider vector bundles $\F$ on $\PP^1\times C$ 
as in the following Lemma.

\begin{lemma}\label{corresp}
Let $f:\PP^1\ra M$ be a morphism and let $\F=(f\times\id_C)^*\U$.
If the rational curve $f$ has degree $k$, then the Chern classes
$c_1(\F)$ and $c_2(\F)$ in $\AA^*(\PP^1\times C)$ are such that
\begin{equation*}
c_1(\F)=k\{pt\}\times C+ \PP^1\times c_1(\xi),\quad
\deg(c_2(\F))=k
\end{equation*}
\end{lemma}

\bp Let $m\in\ZZ$ and $\al\in\AA^1(C)$ such that
\begin{equation*}
c_1(\F)=m\{pt\}\times C+ \PP^1\times\al
\end{equation*}
 
Since for $c\in C$, we have $\det(\U_{|M\times\{c\}})\cong\Th$, it follows 
that
\begin{equation*}
\det(\F_{|\PP^1\times\{c\}})\cong O(k) 
\end{equation*}
 
This implies $m=k$. On the other hand, if $p\in\PP^1$, 
we have $\det(\F_p)\cong\xi$. It follows that $\al=c_1(\xi)$ and 
so we have found the formula for $c_1(\F)$.
From \cite{C} (Lemma 3.1), we have 
\begin{equation*}
\deg(c_1(f^*\Th))=\deg(2c_2(\F)-\frac{c_1(\F)^2}{2})
\end{equation*}

As $f$ has degree $k$, hence, $f^*\Th\cong O(k)$, it follows 
that $\deg(c_2(\F))=k$. \ep

\subsection{Notations}

\

\ps

All our schemes, morphisms and products of schemes are over $\CC$, unless
otherwise stated. If $\E$ is a locally free sheaf on a scheme $X$, 
we denote by $\PP(\E)$ the scheme $\Proj(\Sym(\E^*)$. When working with 
line bundles over $\PP^1\times C$, if $\M$ is a line bundle on $\PP^1$ and 
$\N$ is a line bundle on $C$, we denote by $\M\bo\N$ the sheaf 
$p_1^*\M\otimes p_2^*\N$, where $p_1$ and $p_2$ are the two projections.





\

\textbf{Acknowledgements.}
I thank Johan de Jong for many useful suggestions, and my advisor,
Joe Harris, for all his help, support and inspiration.



\section{Review of Brosius' results}


\subsection{The canonical extension of a rank $2$ vector bundle}\label{can_ext}

\

\ps

Let $T$ be an integral Noetherian scheme. Let $\F$ be a rank $2$
vector bundle on $\PP^1\times T$ and let $k$ be the fiber degree
of $\F$ with respect to the second projection $p_2:\PP^1\times
T\ra T$:
\begin{equation*}
k=c_1(\F). (\PP^1\times \{pt\})
\end{equation*}

If $F$ is a fiber of $p_2$, then the restriction of $\F$ to $F$
has the form $O(a)\oplus O(k-a)$, for some integer $a\geq\frac{k}{2}$.
The pair $(a,k-a)$ is called the \emph{fiber type}  of $\F$ on
$F$.

The function $h:T\ra\ZZ$, which associates to a point $t\in T$ the
integer $a$ from the fiber type of $\F$ (on the fiber $F$ above
$t$), is upper semicontinuos. If $a$ is the value of $h$ at the
generic point of $T$, it is said that $\F$ has \emph{generic fiber
type} $(a,k-a)$. If $a>\frac{k}{2}$, then $\F$ is said to be of type $U$
(unequal) and if $a=\frac{k}{2}$, then $\F$ is said to be of type $E$
(equal).

\ps

Let $\F$ be a rank $2$ vector bundle on $\PP^1\times T$ of generic
fiber type $(a,k-a)$. Denote
\begin{equation*}
\F'=({p_2}^*{p_2}_*(\F(-a)))(a)
\end{equation*}

The canonical morphism $p_2^*{p_2}_*(\F(-a))\ra\F(-a)$ induces a
morphism $g:\F'\ra\F$. It is not hard to see that $\F'$ is a
torsion free sheaf (of rank $1$ in case $U$ and rank $2$ in case
$E$) and that the morphism $g$ is injective. Let $\J=\coker(g)$.
The canonical extension of $\F$ is:
\begin{equation*}
0\ra\F'\ra\F\ra\J\ra0
\end{equation*}

\

\textbf{The canonical extension in Case $U$}

\

In case $U$, the sheaf $\J$ is isomorphic to $\I_Z\otimes\M$, where
$\M$ is some line bundle on $\PP^1\times T$ and $\I_Z$ is the ideal
sheaf of a local complete intersection (lci) subscheme $Z$ of
$\PP^1\times T$. The subscheme $Z$ is given by the zeros of the
morphism $g$, and it has codimension at least $2$ in $\PP^1\times
T$.

\ps

Consider now the case when the scheme $T$ is a smooth projective curve $C$.
Let $\F$ be a rank $2$ vector bundle on $\PP^1\times C$ with
generic fiber type $(a,k-a)$, where $a>\frac{k}{2}$. Let $\Lb$ be
the line bundle ${p_2}_*\F(-a)$. The canonical extension of $\F$
has the form
\begin{equation*}
\begin{CD}
0\ra\F'\ra\F\ra\J\ra0\\\hbox{with} \quad\F'\cong
O(a)\bo\Lb,\quad\M\cong O(k-a)\bo\Lb', \quad\J\cong\I_Z\otimes\M
\end{CD}
\end{equation*}
where $Z\subset\PP^1\times C$ is an lci $0$-cycle and $\Lb'$ some
line bundle. If $\F$ is indecomposable, there are uniquely
determined line bundles $\F'$, $\M$, an lci $0$-cycle
$Z\subset\PP^1\times C$ and a point
\begin{equation*}
\xi\in\PP(\Ext^1(\I_Z\otimes\M, \F'))
\end{equation*}

Such an element $\xi$ is called \emph{invertible} if it determines an
extension with the middle term $\F$ a locally free sheaf. The
invertible elements form a dense open set in
$\PP(\Ext^1(\I_Z\otimes\M, \F'))$.

\ps

The extension determined by an invertible element $\xi$ is, up to
multiplication by a scalar, the canonical extension of $\F$ (see
also Lemma \ref{this_is_can_U}). Moreover, an isomorphism $\F_1\cong\F_2$  of
vector bundles of type $U$ on $\PP^1\times C$ induces an
isomorphism on their canonical sequences. The two sequences
determine the same element in $\PP(\Ext^1_{\PP^1\times
C}(\J,\F'))$.

\

\textbf{The canonical extension in Case $E$}

\

Let $\F$ be a rank $2$ vector bundle on $\PP^1\times T$ with
generic fiber type $(a,a)$. There is a unique line bundle $\M$ on
$\PP^1\times T$, and a closed subscheme $Z\subset\PP^1\times T$
such that, if $D$ is the scheme theoretic image $p_2(Z)\subset T$,
and $\I$ is the ideal sheaf of $Z$ in $\PP^1\times D$, we have
$\J\cong\I\otimes\M$. The subscheme $Z$ is clearly not unique.

\ps

Consider now the case when the scheme  $T$ is a smooth projective curve $C$.
Let $\F$ be a rank $2$ vector bundle on $\PP^1\times C$, with
generic fiber type $(a,a)$. Let $\E$ be the rank $2$ vector bundle
${p_2}_*\F(-a)$. The canonical extension of $\F$ has the form
\begin{equation*}
\begin{CD}
0\ra\F'\ra\F\ra\J\ra0,\quad\hbox{with } \quad\F'\cong
O(a)\bo\E,\quad \J\cong\I(a)
\end{CD}
\end{equation*}
where  $\I$ is the ideal sheaf of a $0$-cycle $Z$ in $\PP^1\times
D$ -- where we denote by $D$  the scheme theoretic image $p_2(Z)$.

\ps

We restrict ourselves to analyze canonical extensions of vector
bundles $\F$ such that $\E$ is a \emph{stable} vector bundle on
$C$ and the morphism $p_2:Z\ra D$ is an isomorphism, in which
case we say that $\F$  has type $(\da)$. Note that we have:
\begin{equation}\label{type_da}
\I\cong O(-1)\bo O_D,\quad\quad\J\cong O(a-1)\bo O_D
\end{equation}

If $\F$ is indecomposable, there are uniquely determined sheaves
$\F'$, $\J$ and an orbit $\xi$ for the action of
$\Aut(O_D)\cong\H^0(D,O^*_D)$ on $\Ext^1(\J, \F')$.

An orbit $\xi$ is called \emph{invertible} if it determines an
extension with the middle term $\F$ a locally free sheaf. The
invertible elements form a dense open set in $\Ext^1(\J, \F')$.

\ps

The extension determined by an invertible element $\xi$ and the canonical
extension of its middle term are in the same orbit of the action
of $\Aut(O_D)$ on $\Ext^1(\J, \F')$. Moreover, if there is an
isomorphism $\F_1\cong\F_2$  of vector bundles of type $(\da)$ on
$\PP^1\times C$, then their canonical extensions are in the same
orbit for the action of $\Aut(O_D)$ on $\Ext^1(\J, \F')$.

\

\ps

One may construct moduli spaces for rank $2$ vector bundles $\F$
on $\PP^1\times C$ using the canonical extension of $\F$. Let
$\xi$ be a fixed line bundle of degree $1$. Let $k\geq1$ be an
integer and denote:
\begin{equation}\label{c}
c_1=k\{pt\}\times C+ \PP^1\times c_1(\xi),\quad\quad c_2=k
\end{equation}


\subsection{The moduli space in Case $U$}

\

\ps

Let $a$ and $e$ be integers such that $a>\frac{k}{2}, e\geq0$. We
say that a line bundle on $\PP^1\times C$ has \emph{type} $(a,-e)$
if it has the form $O(a)\bo\Lb$, where $\Lb$ is a line bundle on
$C$ of degree $-e$. Define a contravariant functor
\begin{equation}\label{F_odd_case}
\begin{CD}
F:\Sch_{\CC}\lra\Sets\\
F(S)=\{\hbox{Isomorphism classes of rank $2$ vector bundles $\F$
on $\PP^1\times C\times S$, such that}\\
\hbox{$\forall s\in S$, the bundle $\F_s$ has Chern classes $c_1$ and
$c_2$ as in (\ref{c})}\\
\hbox{and the canonical subbundle of $\F_s$ has type $(a,-e)$}\}
\end{CD}
\end{equation}

Note that if $\F\in F(\Spec(\CC))$, then the canonical extension of $\F$ has
the form:
\begin{equation}\label{U_ext}
0\ra
O(a)\bo\Lb\ra\F\ra\I_Z\otimes(O(k-a)\bo(\Lb^{-1}\otimes\xi))\ra0
\end{equation}
where $\Lb\in\Pic^{-e}(C)$ and $\I_Z$ is the ideal sheaf of a
$0$-cycle $Z$. If we let $\de=\length(Z)$, then, by a computation with Chern
classes, we have:
\begin{equation}\label{de}
\de=(k-a)-e(2a-k)
\end{equation}

Note that $e\geq0$ and $\de\geq0$ imply $\frac{k-a}{2a-k}\geq e\geq0$ and
$a\leq k$.

\begin{thm}\cite{BR2}\label{B(a,e)}
For each pair of integers $a$ and $e$ such that $k\geq
a>\frac{k}{2}, e\geq0$, there exists an irreducible variety $\Bb$
which is a fine moduli scheme for the functor $F$ of
(\ref{F_odd_case}).
\end{thm}

Precisely, there is a universal rank $2$ bundle $\F_{\Bb}$ on
$\PP^1\times C\times\Bb$: if $S$ is a scheme  and $\F_S\in F(S)$,
then there is a unique morphism $\nu:S\ra\Bb$ such that
$\F_S\cong\nu^*\F_{\Bb}\otimes\N$, for some line bundle $\N$ on
$S$.

\

\textbf{Outline of the construction of $\Bb$}

\

There is an irreducible  projective variety $\Bb'$ parameterizing
extensions as in (\ref{U_ext}). If we denote by
$\Hilb^{\de}(\PP^1\times C)$ the Hilbert scheme of lci
$0$-dimensional subschemes of $\PP^1\times C$ having length $\de$,
$\Bb'$ is a projective bundle over the scheme
$\Pic^{-e}(C)\times\Hilb^{\de}(\PP^1\times C)$, which we denote with $\Ss$.
Let $p':\Bb'\ra\Ss$ be the canonical map. The fibers of $p'$ are of
the form $\PP(\Ext^1(\I_Z\otimes\M,\F'))$, where $\F'$ and $\M$
are from the canonical extension (\ref{U_ext}) of $\F$.

\ps

There is a universal bundle $\F_{\Bb'}$ on $\PP^1\times
C\times\Bb'$ and a dense  open $\Bb\subset\Bb'$, such that
$\F_{\Bb'}$ restricts to a locally free sheaf $\F_{\Bb}$ on $\Bb$.
Given the observations in \ref{can_ext}, one may prove that
$\Bb$ is a fine moduli space for the functor $F$ of
(\ref{F_odd_case}). Note that $\Bb$ intersects every fiber of
$p':\Bb\ra\Ss$ in a dense open set.

\ps

\begin{notn}
For integers $a$ and $e$ as in Theorem \ref{B(a,e)} we denote by
$\Bb(a,e)$ the moduli scheme $\Bb$; we denote by $\Bb'(a,e)$ the
variety $\Bb'$.
\end{notn}

The variety $\Bb'(a,e)$ is a $\PP^N$ bundle over
$\Ss=\Pic^{-e}(C)\times\Hilb^{\de}(\PP^1\times C)$, where $\de$ is
given by (\ref{de}) and $N$ is given by
\begin{equation*}
N=\dim\Ext^1(\I_Z\otimes\M,\F')
\end{equation*}

Note that we have $N=\chi(\F'\otimes\M^*)+\de-1$. In
terms of $a$, $e$ and $k$, this is:
\begin{equation*}
N=(2a-k+1)(2e+g)+\de=(2a-k+1)g+(2a-k+2)e+(k-a)
\end{equation*}

It follows that the dimension of $\Bb(a,e)$ is
\begin{equation}\label{dim_B(a,e)}
\dim\Bb(a,e)=2\de+g+N-1=(2a-k+2)g+(3k-3a-1)-e(2a-k-2)
\end{equation}

Note that $N>0$ and the variety $\Bb(a,e)$ is non-empty for every
$a$ and $e$ in the given range.
\begin{rmk}\label{extra_B}
If $k$ is even, one may allow  $a=\frac{k}{2}$ in the construction
of the variety $\Bb(a,e)$. In this case, $\de=\frac{k}{2}$ and
$\Bb(\frac{k}{2},e)$ is non-empty for any $e\geq0$. The
variety $\Bb(\frac{k}{2},e)$ is not a fine moduli scheme for the
functor $F$.
\end{rmk}


\subsection{The moduli space in Case $E$}

\

\ps

Let $k$ be a positive even integer and let $a=\frac{k}{2}$.
Define a contravariant functor
\begin{equation}\label{F_even_case}
\begin{CD}
F:\Sch_{\CC}\lra\Sets\\
F(S)=\{\hbox{Isomorphism classes of rank $2$ vector bundles $\F$
on $\PP^1\times C\times S$, such that}\\
\hbox{$\forall s\in S$, the bundle $\F_s$ has Chern classes $c_1$
and $c_2$ as in (\ref{c}) and type $(\da)$}\}
\end{CD}
\end{equation}

\

If $\F\in F(\Spec(\CC))$, then the canonical extension of $\F$ has
the form:
\begin{equation}\label{E_ext}
0\ra O(a)\bo\E\ra\F\ra O(a-1)\bo O_D\ra0
\end{equation}
where $D$ is a $0$-cycle on $C$ of length $a$ and $\E$ is a
stable, rank $2$ vector bundle on $C$, with determinant
$\det(\E)\cong\xi(-D)$.
\begin{thm}\cite{BR2}\label{B_E}
For every positive even integer $k$, there exists an irreducible
variety $\Bb$ which is a coarse moduli scheme for the functor $F$
of (\ref{F_even_case}).
\end{thm}

It follows that there is a functorial way of associating to a scheme $S$ and
an element $\F_S\in F(S)$, a morphism $\nu:S\ra\Bb$.

\

\textbf{Outline of the construction of $\Bb$}

\

Let $M_{1-a}$ be the coarse moduli scheme of rank $2$, semistable
vector bundles on $C$, with determinant of degree $(1-a)$. Let the
determinant map be $det:M_{1-a}\ra\Pic^{1-a}(C)$. Let
$M^s_{1-a}\subset M_{1-a}$ be the locus of stable vector bundles.

\ps

The moduli space $M_{1-a}$ is the geometric quotient of a smooth,
quasiprojective variety $\MM_{1-a}$, by the action of an algebraic
group $\PGL(r)$, for some integer $r$. Let $\tau:\MM_{1-a}\ra
M_{1-a}$ be the quotient map and let
$\MM^s_{1-a}=\tau^{-1}(\MM^s_{1-a})$. Note that
there is no {Poincar\'e} bundle on $M_{1-a}\times C$, not even on an open
set \cite{R}. However, there is a {Poincar\'e} bundle $\V$ on
$\MM_{1-a}\times C$.

\ps

Define $\th:\Pic^{1-a}(C)\ra\Pic^a(C)$ to be the morphism that
sends $\Lb$ to $\xi\otimes\Lb^{-1}$. Consider the composition
$\th\circ(det)\circ\tau:\MM^s_{1-a}\ra\Pic^a(C)$. We define
$\SSs$ to be the fiber product
\begin{equation*}
\SSs=\MM^s_{1-a}\times_{\Pic^a(C)}\Sym^a(C)
\end{equation*}

\ps

There is a vector bundle $\BBb'$ over $\SSs$, parametrizing
extensions as in (\ref{E_ext}). Denote by $\pp':\BBb'\ra\SSs$  the
canonical map. The fibers of $\pp'$ are of the form
$\Ext^1(\J,\F')$, where $\F'$ and $\J$ are from the canonical
sequence (\ref{E_ext}) of $\F$.

\ps

There is a locally universal bundle $\F_{\BBb'}$ on $\PP^1\times
C\times\BBb'$: if $S$ is a scheme  and $\F_S\in F(S)$, there is a
covering of $S$ with open sets $S_i$, and morphisms
$\nu_i:S_i\ra\BBb$  such that
$\F_{S_i}\cong\nu_i^*\F_{\BBb}\otimes\N$, for some line bundle
$\N$ on $S_i$. Moreover, there is a dense  open
$\BBb\subset\BBb'$, such that $\F_{\BBb'}$ restricts to a locally
free sheaf $\F_{\BBb}$ on $\BBb$. Note that $\BBb$ intersects
every fiber of $\pp':\BBb'\ra\Ss$ in a dense open set.

\ps

Let $T$ be the group scheme over $\SSs$ which has as fibers over
$\SSs$ the algebraic groups given by $\Aut(O_D)$. Then $T$ acts on
$\BBb$ by the action of $\Aut(O_D)$ on $\Ext^1(\J,\F')$. If $G'$
is the group $\PGL(r)$ acting on $\MM_{1-a}$, then $G'$ acts as well on
$\BBb$. Identifying the scalar multiplication of extensions that appears
in both the action of $T$ and $G'$, one may define a group scheme $G$
over $\SSs$, which acts on $\BBb$. There is a universal geometric quotient
$\Bb$ for the action of $G$ on $\BBb$, which make the scheme $\Bb$
into a  coarse moduli scheme for the
functor $F$ of (\ref{F_even_case}).

There is a geometric quotient $\Ss$ of $\SSs$ by $G$, which is
\begin{equation*}
\Ss\cong M^s_{1-a}\times_{\Pic^a(C)}\Sym^a(C)
\end{equation*}

The projection map $\pp'$ descends to give a map $p:\Bb\ra\Ss$.
Away from the ramification locus of $C^a\ra\Sym^a(C)$, the fibers of the
map $p$ are isomorphic to $({\PP}GL(1))^a$.

\begin{notn}
We denote by $\Bbe$ the moduli scheme $\Bb$ from Theorem \ref{B_E}.
\end{notn}

The dimension of the variety $\Ss$ is $(a+3g-3)$. It follows that
\begin{equation}
\dim(\Bbe)=4a+3g-3=2k+3g-3.
\end{equation}


\section{The good locus in the moduli space of bundles on $\PP^1\times C$}

\subsection{The good locus in case $U$}

\

\ps

Let $\Bb=\Bb(a,e)$ be the variety in Theorem \ref{B(a,e)} and let
$\F_{\Bb}$ be the universal bundle on $\PP^1\times C\times\Bb$. If
$b\in\Bb$, denote $\F_b=\F_{|\PP^1\times C\times\{b\}}$
and $\F_{p,b}={\F_b}_{|\{p\}\times C}$.

We say that a point $b\in\Bb$ is \emph{good} if the bundle $\F_b$
induces a morphism $f:\PP^1\ra M$, or, equivalently, the bundle $\F_{p,b}$ is
stable for any $p\in\PP^1$. Let $\Bb^0\subset\Bb$ be the set of good points.
The property of being stable is an open condition, so $\Bb^0$ is
open in $\Bb$, but possibly empty.

\begin{thm}\label{good_B(a,e)}
Let $a$ and $e$ be integers such that $a\geq\frac{k}{2}$ and $\de\geq0$,
where $\de=(k-a)-e(2a-k)$. The open $\Bb^0(a,e)$ is non-empty if and
only if either $e>0$ or if $e=0$ and $\de=0$
\end{thm}

\bp
We first prove that the conditions are necessary. Consider an element
in $\Bb(a,e)$:
\begin{equation*}
0\ra O(a)\bo\Lb\ra\F\ra (O(k-a)\bo(\Lb^{-1}\otimes\xi))\otimes\I_Z\ra0
\end{equation*}
where $\Lb\in\Pic^{-e}(C)$ and $Z\subset\PP^1\times C$ is a $0$-cycle
of length $\de$.

If $p\in\PP^1$, then it follows that there is an injective morphism
$\Lb\ra\F_p$. As the bundle  $\F_p$ has degree $1$, if it is stable,
then we must have $-e\leq0$. It follows that $\Bb^0(a,e)=\eset$ if $e<0$.

\ps

Assume now that $e=0$ and $\de>0$. Let $p\in\PP^1$ be such that
$(\{p\}\times C)\cap Z\neq\eset$ and let $D$ be the $0$-cycle on $C$ such that
$(\{p\}\times C)\cap Z=\{p\}\times D$. By restriction to $\{p\}\times C$
we get:
\begin{equation*}
0\ra\Lb\ra\F_p\ra(\Lb^{-1}\otimes\xi(-D))\oplus(\Lb^{-1}\otimes\xi_{|D})\ra0
\end{equation*}

By passing to the saturation, we have that $\Lb(D)$ is a subbundle of $\F_p$.
As $\deg(\Lb(D))=\deg(D)>0$, it follows that $\F_p$ is not stable, thus showing
that $\Bb^0(a,e)=\eset$.

\ps

We prove now that the conditions are sufficient. From the construction of
$\Bb=\Bb(a,e)$, recall the notation
\begin{equation*}
\Ss=\Pic^{-e}(C)\times\Hilb^{\de}(\PP^1\times C)
\end{equation*}
and consider the projective bundle $p':\Bb'\ra\Ss$.

We prove that $\Bb^0$ is dense in a general fiber of $p'$.
Fix $\si=(\Lb,Z)\in\Ss$. Denote
\begin{equation*}
\F'=O(a)\bo\Lb, \quad\M=O(k-a)\bo(\Lb^{-1}\otimes\xi)
\end{equation*}

The fiber  ${p'}^{-1}(\si)$ is isomorphic to a projective space
$\PP(W)$, where we denote
\begin{equation*}
W=\Ext^1_{\PP^1\times C}(\M\otimes\I_Z,\F')
\end{equation*}

Denote with $U$ the set $\Bb\cap {p'}^{-1}(\si)$, which is the dense open
set in $\PP(W)$ corresponding to extensions which have the middle term a locally
free sheaf.  Let $U^0\subset U$ be the good locus.

Define
\begin{equation*}
\begin{CD}
Y=\{(p,u)\in\PP^1\times U|
\hbox{ the bundle } \F_{p,u} \hbox{ is not stable }\}\\
Y_p=\{u\in U| \hbox{ the bundle } \F_{p,u} \hbox{ is not stable }\}
\end{CD}
\end{equation*}

The set $Y$ is closed in $\PP^1\times U$ and, for any $p\in\PP^1$, the set
$Y_p$ is closed in $U$. If $\pi:\PP^1\times U\ra U$ is the second projection,
we have $U^0=U\setm\pi(Y)$.

\ps

We prove the following: \bi \item[i) ]for any $u\in\Ss$, the open
$U^0$ is not empty if $\de=0$ \item[ii) ]for general $u\in\Ss$,
the open $U^0$ is not empty if $e>0$ and $\de>0$
\ei

In both cases, we show that $Y$ has codimension at least $2$ in
$\PP^1\times U$, therefore, $\pi(Y)$ is a proper closed subset of
$U$ and $U^0\neq\eset$.

\ps

\emph{Proof of i).} We prove that for any $p\in\PP^1$, the set
$Y_p$ has codimension at least $2$ in $U$.

\ps

Since $\de=0$ we have $Z=\eset$ and $W=\Ext^1_{\PP^1\times
C}(\M,\F')$. Denote
\begin{equation*}
 V=\Ext^1_{\{p\}\times C}(\M_{|\{p\}\times C}, \F'_{|\{p\}\times C}).
\end{equation*}

Let $r:W\ra V$ be the restriction morphism, which, by  Lemma
\ref{restr_1}, is surjective. Note
\begin{equation*}
\M_{|\{p\}\times C}\cong\Lb^{-1}\otimes\xi,\quad \F'_{|\{p\}\times
C}\cong\Lb
\end{equation*}

It follows that there is an isomorphism
$V\cong\Ext^1(\Lb^{-1}\otimes\xi,\Lb)$. From \cite{C} (Lemma 2.1)
we have that the locus of unstable extensions in $\PP(V)$
is a closed subvariety $Z\subset\PP(V)$  of codimension at least
$g\geq2$. Let $\tilde{Z}\subset\PP(W)$ be the closure of the
preimage of $Z$ via the rational map $\PP(W)\dra\PP(V)$ induced by
the restriction map $r:W\ra V$. Then $Y_p=\tilde{Z}\cap U$ and we
have $\codim_U Y_p=\codim_{\PP(V)}Z\geq g\geq2$.

\ps

\emph{Proof of ii).} Since $\de>0$, we have $Z\neq\eset$. Let
$u=(\Lb,Z)\in\Ss$ be a general point such that $Z$ is reduced,
with no two points in the same fiber of the first projection
$\PP^1\times C\ra\PP^1$. Let $\Ga\subset\PP^1$ be the image of $Z$
in $\PP^1$. Note that if $p\in\PP^1\setm\Ga$, we have
${\I_Z}_{|\{p\}\times C}\cong O$. The same argument as in i)
proves that $\codim_U Y_p\geq2$.

\ps

We prove that $Y_p$ has codimension at least $1$ if $p\in\Ga$.
Equivalently, a general extension in $W=\Ext^1_{\PP^1\times
C}(\M\otimes\I_Z,\F')$ has its middle term $\F$ so that $\F_p$ is
a stable bundle. Assume the contrary: any extension in $W$ has
middle term $\F$ such that $\F_p$ is unstable. Denote
\begin{equation*}
V=\Ext^1_{\{p\}\times C}((\M\otimes\I_Z)_{|\{p\}\times
C},\F'_{|\{p\}\times C})
\end{equation*}

Consider the restriction morphism $r:W\ra V$; by Lemma
\ref{restr_1}, the map $r$ is surjective. It follows that any
extension in $V$ has middle term an unstable vector bundle. We
show that this is impossible. Let $q\in C$ be such that
$\{(p,q)\}=(\{p\}\times C)\cap Z$. Note
\begin{equation*}
\F'_{|\{p\}\times C}\cong\Lb,\quad(\M\otimes\I_Z)_{|\{p\}\times
C}\cong\Lb^{-1}\otimes\xi\otimes(O(-q)\oplus O_q)
\end{equation*}

Denote
\begin{equation*}
\Lb'=\Lb(q), \quad V'=\Ext^1_C({\Lb'}^{-1}\otimes\xi,\Lb')
\end{equation*}

Since $e>0$, we have that $\deg(\Lb')\leq0$ and we can apply again
Lemma 2.1 from \cite{C}, to get that a general extension in $V'$
has the middle term a stable vector bundle.

Consider an exact sequence
\begin{equation}\label{ext'}
0\ra\Lb'\ra\E\ra{\Lb'}^{-1}\otimes\xi\ra0
\end{equation}

By composing $\Lb'\ra\E$ with the inclusion $\Lb\ra\Lb'$, we
obtain an exact sequence
\begin{equation}\label{ext}
0\ra\Lb\ra\E\ra(\Lb^{-1}\otimes\xi\otimes(O(-q))\oplus O_q\ra0
\end{equation}

But if the extension (\ref{ext'}) is general in $V'$, we have that
$\E$ is stable. But the extension (\ref{ext}) is in $V$ and it has
middle term the bundle $\E$, contradiction. This completes the
proof. \ep

\

We give now the proof of the auxiliary lemma used in the proof of
Theorem \ref{good_B(a,e)}.

\begin{lemma}\label{restr_1}
Let $X$ be an irreducible variety and $Y\subset X$ a closed
irreducible subvariety with ideal sheaf $\I$.
Let $\G'$ and $\G''$ coherent sheaves on $X$ such that
for any exact sequence
\begin{equation*}
0\ra\G'\ra\G\ra\G''\ra0
\end{equation*}
by restriction to $Y$ we get an exact sequence of coherent sheaves on $Y$:
\begin{equation*}
0\ra\G'\otimes O_Y\ra\G\otimes O_Y\ra\G''\otimes O_Y\ra0
\end{equation*}
The restriction gives a linear map:
\begin{equation*}
r:\Ext^1_X(\G'',\G')\ra\Ext^1_Y(\G''\otimes O_Y,\G'\otimes O_Y)
\end{equation*}
Moreover, if $\Ext^2_X(\G'',\G'\otimes\I)=0$,
the morphism $r$ is surjective.
\end{lemma}

The conditions in Lemma \ref{restr_1} are satisfied if, for
example, the sheaf $\G'$ is locally free and the sheaf $\G''$ has
the property that $\rk\G''=\rk_Y(\G''\otimes O_Y)$. (If $X$ is an
irreducible variety, $\eta$ is its generic point and $\G$ is a
coherent sheaf on $X$, we define the rank $\rk\G$ as the dimension
of the $O_{\eta}$-vector space $\G_{\eta}$).

\bp It is clear that the map $r$ commutes with scalar
multiplication. It is straightforward to check that the map $r$
preserves addition, using the definition of the sum of two
extensions (see \cite{W}, Def. 3.4.4). Consider the canonical
morphisms of vector spaces:
\begin{equation*}
\begin{CD}
w:\Ext^1_X(\G'',\G')\ra\Ext^1_Y(G'',\G'\otimes O_Y)\\
v:\Ext^1_Y(\G''\otimes O_Y,\G'\otimes O_Y)
\ra\Ext^1_Y(G'',\G'\otimes O_Y)
\end{CD}
\end{equation*}

It is easy to check that $v\circ r=w$. Since we have
$\coker(w)\cong\Ext^2_X(\F'',\F'\otimes\I)$, the lemma follows.
\ep


\subsection{The good locus in case $E$}

\

\ps

Let $\Bb=\Bbe$ be the variety in Theorem \ref{B_E}. Recall from the
construction of $\Bb$, that it is the geometric quotient of an irreducible 
variety  $\BBb$ by the action of a group $G$. Let $\F_{\BBb}$ be the 
locally universal bundle on $\PP^1\times C\times\BBb$. If $b\in\BBb$, 
denote by $\F_b=\F_{|\PP^1\times C\times\{b\}}$ and
$\F_{p,b}={\F_b}_{|\{p\}\times C}$.

We say that a point $b\in\BBb$ is \emph{good} if the bundle $\F_b$
induces a morphism $f:\PP^1\ra M$, or, equivalently, the bundle
$\F_{p,b}$ is stable for any $p\in\PP^1$. Let $\BBb^0\subset\BBb$
be the set of good points. The property of being stable is an open
condition, so $\BBb^0$ is open in $\BBb$, but possibly empty.
\begin{thm}\label{good_B_E}
The open $\BBb$ is non-empty.
\end{thm}

\bp
This is similar to the proof of Theorem \ref{good_B(a,e)}. Let $k=2a$.
From the construction of $\BBb$, recall the notation
\begin{equation*}
\SSs=\MM^s_{1-a}\times_{\Pic^a(C)}\Sym^a(C)
\end{equation*}
and consider the vector bundle $\pp':\BBb'\ra\SSs$.

We prove that $\BBb^0$ is dense in a general fiber of $\pp'$. Let
$\si\in\Ss$ and let $(\E,D)$ be its image (by the canonical
projection) in $M^s_{1-a}\times_{\Pic^a(C)}\Sym^a(C)$. Let
\begin{equation*}
\F'=O(a)\bo\E, \quad\quad\M=O(a-1)\bo O_D
\end{equation*}

The fiber  ${\pp'}^{-1}(\si)$ is isomorphic to the affine space
$W=\Ext^1_{\PP^1\times C}(\J,\F')$. Denote with $U$ the
set $\BBb\cap {\pp'}^{-1}(\si)$, which is the dense open in $W$ corresponding
to extensions whose middle term is a locally free sheaf. Let
$U^0\subset U$ be the good locus.

Define:
\begin{equation*}
\begin{CD}
Y=\{(p,\si)\in\PP^1\times U|
\hbox{ the bundle }  \F_{p,\si} \hbox{ is not stable }\}\\
Y_p=\{\si\in U| \hbox{ the bundle } \F_{p,\si} \hbox{ is not stable }\}
\end{CD}
\end{equation*}

The set $Y$ is closed in $\PP^1\times U$ and, for any
$p\in\PP^1$, the set $Y_p$ is closed in $U$. Let $\pi:\PP^1\times
U\ra U$ be the second projection. Then $U^0=U\setm\pi(Y)$.

\ps

We prove that if $\si\in\SSs$ is a general element such that $D\in\Sym^a(C)$
consists of distinct points, then $Y$ has codimension at least $2$
in $\PP^1\times U$, and therefore, $\pi(Y)$ is a proper closed
subset of $U$. We prove this by showing that for any $p\in\PP^1$,
the set $Y_p$ has codimension at least $2$ in $U$.

\ps Let $p\in\PP^1$ and denote
\begin{equation*}
 V=\Ext^1_{\{p\}\times C}(\J_{|\{p\}\times C}, \F'_{|\{p\}\times C}).
\end{equation*}

Let $r:W\ra V$ be the restriction morphism, which, by  Lemma
\ref{restr_1}, is surjective. Note
\begin{equation*}
\F'_{|\{p\}\times C}\cong\E,\quad \J_{|\{p\}\times C}\cong O_D
\end{equation*}

From \cite{C} (Lemma 5.4) we have that the locus of
unstable extensions in $\PP(V)$ is a closed subvariety
$Z\subset\PP(V)$  of codimension at least $2$. Let
$\tilde{Z}\subset W$ be the preimage of the affine cone over $Z$
via the restriction morphism $r$. Then $Y_p=\tilde{Z}\cap U$ and
we have:
\begin{equation*}
\codim_U Y_p=\codim_{\PP(V)}Z\geq2
\end{equation*}
\ep

\begin{rmk_def}
The action of the group $G$ on $\BBb$ induces an action
on $\BBb^0$. Define the good locus $\Bb^0\subset\Bb$
as the image of the good locus $\BBb^0$ via the quotient map.
It follows from Theorem \ref{good_B_E} that $\Bb^0\neq\eset$.
\end{rmk_def}


\section{The loci $\M(a,e)$ and $\M_E$ in the space of rational curves}

We find all the irreducible components of the space
$\Hom_k(\PP^1,M)$. They all correspond to the moduli of rank $2$
vector bundles on $\PP^1\times C$ that we constructed.

\ps

Let $k\geq1$ be an integer and let $\Bb$ be one of the moduli spaces
$\Bbe$ or $\Bb(a,e)$, for $a$ and $e$ such that the good locus is
non-empty (we assume $k$ even if we work with $\Bbe$).
Let $F$ be the functor for which $\Bb$ is a coarse moduli scheme.
Define $F^0$ to be the following open subfunctor of $F$:
\begin{equation}\label{F^0}
\begin{CD}
F^0:\Sch_{\CC}\lra\Sets\\
F^0(S)=\{\F\in F(S)|
\quad\forall\mbox{ } (p,s)\in\PP^1\times S, \hbox{ the bundle }\F_{p,s}
\hbox{ is stable }\}
\end{CD}
\end{equation}

The good locus $\Bb^0\subset\Bb$ is a coarse moduli scheme for the
functor $F^0$.

\ps

Consider the evaluation (universal) morphism
\begin{equation*}
\ev:\PP^1\times\Hom_k(\PP^1,M)\ra M
\end{equation*}

An element $\G\in F(S)$ induces a morphism $g:\PP^1\times S\ra M$.
There is a unique morphism $h:S\ra\Hom_k(\PP^1,M)$ such that
$g=\ev\circ(\id_{\PP^1}\times h)$. There is a transformation of functors:
\begin{equation}\label{T}
T:F^0\lra\Hom(-,\Hom_k(\PP^1,M))
\end{equation}
which associates to an element $\G\in F^0(S)$ the element
$h\in\Hom(S,\Hom_k(\PP^1,M))$.

\begin{rmk}\label{T'}
There is a morphism
\begin{equation}\label{rho}
\rho:\Bb^0\ra\Hom_k(\PP^1,M)
\end{equation}
that sends an element an element in $\Bb^0$ that corresponds to a vector bundle
$\F$ on $\PP^1\times C$, to the morphism $f:\PP^1\ra M$ induced by $\F$. The
transformation $T$ factors through
\begin{equation*}
T':F^0\lra\Hom(-,\Bb^0)
\end{equation*}
This is an immediate consequence of the coarse moduli property of $\Bb$.
\end{rmk}

\begin{lemma}
The morphism $\rho$ is birational onto its image.
\end{lemma}

\bp
This follows immediately from the fact that a
morphism $f:\PP^1\ra M$ is induced by a unique $\F$ on
$\PP^1\times C$; hence, $\rho$ is an injective map.
\ep

\

If $\U$ is the rigidified Poincar\'e bundle on $M\times C$, let $\HH$ be the
vector bundle $(ev\times id_C)^*\U$
on $\PP^1\times C\times\Hom_k(\PP^1,M)$.
The bundle $\HH$ is a universal bundle: if $f\in\Hom_k(\PP^1,M)$,
then we have
\begin{equation*}
\HH_f=\HH_{|\PP^1\times C\times \{f\}}\cong(f\times{id})^*\U
\end{equation*}

\ps

If $\M\subset\Hom_k(\PP^1,M)$ is an irreducible component (with
the reduced structure), we denote by $\HH_{\M}$ the restriction of
the universal bundle $\HH$ to $\PP^1\times C\times\M$.

\begin{thm}\label{M^0}
Let $\M\subset\Hom_k(\PP^1,M)$ be an irreducible component. There
is a dense open $\M^0\subset\M$ and a moduli space $\Bb$
such that the members of the family $\HH_{\M^0}$, of vector bundles on
$\PP^1\times C$, correspond to points in $\Bb^0\subset\Bb$.
\end{thm}

Equivalently, Theorem \ref{M^0} says that if $F^0$ is the functor in
(\ref{F^0}) corresponding to $\Bb$, then we have $\HH_{\M^0}\in F^0(\M^0)$.
Clearly, the moduli space $\Bb$, associated to $\M$,  is
unique. We prove Theorem \ref{M^0} in Section \ref{pf_M^0} and give here
some consequences.

\begin{cor}\label{M_B_corresp}
Let $\M\subset\Hom_k(\PP^1,M)$ be an irreducible component.
Consider the open set $\M^0\subset\M$ and the moduli scheme $\Bb$
corresponding to it, as in Theorem \ref{M^0}. Then the morphism
$\rho$ of (\ref{rho}) gives an isomorphism of the good locus $\Bb^0$
to $\M^0$. This gives a one-to-one correspondence between irreducible
components of $\Hom_k(\PP^1,M)$ and moduli spaces of vector bundles on
$\PP^1\times C$ for which the good locus is not empty.
\end{cor}

\bp
Consider the following map from the transformation (\ref{T}):
\begin{equation*}
T(\M^0):F^0(\M^0)\lra\Hom(\M^0,\Hom_k(\PP^1,M))
\end{equation*}

The element $\HH_{\M^0}\in F^0(\M^0)$ corresponds by the map $T(\M^0)$ to
the inclusion morphism $i:\M^0\ra\Hom_k(\PP^1,M)$. By Remark \ref{T'},
it follows that there is a morphism $\nu:\M^0\ra\Bb^0$ such that
$\rho\circ\nu=i$.
Therefore, we have $\M^0\cong\Bb^0$.
\ep

\begin{defn}
Let $\M(a,e)$ (resp. $\M_E$) be the irreducible
variety which is the image of the morphism $\rho$ in Remark
\ref{T'} for $\Bb=\Bb(a,e)$ (resp. $\Bb=\Bbe$).
Both varieties are considered with the reduced induced structure.
\end{defn}

\begin{cor}\label{M_isom_B^0}
There  are isomorphisms $\M(a,e)\cong\Bb^0(a,e)$ and
$M_E\cong\Bbe^0$ and we have:
\begin{equation}\label{dim_M(a,e)}
\dim\M(a,e)=(2a-k+2)g+(3k-3a-1)-e(2a-k-2)
\end{equation}
\begin{equation}\label{dim_M_E}
\dim\M_E=2k+3g-3
\end{equation}
\end{cor}

The subvarieties $\M(a,e)$ (for all $a$ and $e$) and $\M_E$ (if
$k$ is even) are mutually disjoint and they cover
$\Hom_k(\PP^1,M)$.

\begin{cor}
Any irreducible component of $\Hom_k(\PP^1,M)$ is either $\ME$ or
$\Mae$, for some $a$ and $e$ for which the good locus is not empty.
\end{cor}

\begin{prop_def}\label{M_nice}(\textbf{The nice component})
For any $k\geq1$ there is an irreducible component $\M$ of
$\Hom_k(\PP^1,M)$, of the expected dimension $(2k+3g-3)$, given by:
\bi \item[i) ]If $k=2a+1$ then
$\M$ is the closure in $\Hom_k(\PP^1,M)$ of the variety $\M(a+1,a)$
\item[ii) ]If $k=2a$ then $\M$ is the closure in $\Hom_k(\PP^1,M)$ of the
variety $\M_E$
\ei

We call the component $\M$ the \emph{nice component} of the space of rational
curves $\Hom_k(\PP^1,M)$.
\end{prop_def}

\bp If $k=2a+1$, for $f\in\M(a+1,a)$ the bundle $\HH_f$ has fiber
type $(a+1,a)$. If $k=2a$, for $f\in\M_E$, the bundle $\HH_f$ has
fiber type $(a,a)$. By \cite{C} (Lemma 3.4), $f$ is an
unobstructed point of $\Hom_k(\PP^1,M)$. Therefore, there is a
unique irreducible component $\M$ containing $f$ and, moreover,
$\M$ has the expected dimension $(2k+3g-3)$. Since by
(\ref{dim_M(a,e)}) and (\ref{dim_M_E}) both $\M(a+1,a)$ and $\ME$
have the expected dimension, our result follows.
\ep

\begin{cor}\label{nice_comp_even}
If $k=2a$ is  an even integer, the variety $\ME$ is the unique
irreducible component with general element $f$ such that $\HH_f$
has generic fiber type $(a,a)$.
\end{cor}

\bp From Corollary \ref{M_B_corresp}, it follows that any
irreducible component $\M\subset\Hom_k(\PP^1,M)$ with the property
that a general element $f\in\M$ is such that $\HH_f$ has generic
fiber type $(a,a)$, has an open set $\M^0$ which is isomorphic to
$\Bbe^0$ and the composition
\begin{equation*}
\begin{CD}
\M^0 @>{\nu}>>\Bbe^0 @>{\rho}>>\Hom_k(\PP^1,M)
\end{CD}
\end{equation*}
is the inclusion morphism. It follows that $\M^0=\M_E$; hence, 
there can not be two such components $\M$. By Proposition \ref{M_nice}, the 
subvariety $\ME$ is an irreducible component.
 \ep

From the constructions of the varieties $\Bb(a,e)$ and $\M(a,e)$, 
we have a similar remark for the odd degree case.

\begin{rmk}\label{nice_comp_odd}
If $k=(2a+1)$ is an odd integer, the variety $\Mae$ is the unique
irreducible component with general element $f$ such that $\HH_f$
has generic fiber type $(a+1,a)$.
\end{rmk}

Note that by Lemma \ref{nice_comp_even} and Remark \ref{nice_comp_odd},
it follows that the nice component, as defined in Proposition \ref{M_nice},
is the same as the nice component defined in  \cite{C}. See
Lemma \ref{geom_descr} for another more geometric proof of this fact.


\section{Families of vector bundles on $\PP^1\times C$}

We analyze the behavior of the canonical extension in families of
vector bundles on  $\PP^1\times C$. In this section $S$ will be an
arbitrary irreducible variety. If $\F$ is a vector bundle on
$\PP^1\times C\times S$, we denote $\F_s=\F_{|\PP^1\times
C\times\{s\}}$ and $\F_{c,s}=\F_{|\{c\}\times\{s\}}$, for $s\in
S$, $c\in C$. We denote by $p^s_2$ the canonical projection
$\PP^1\times C\times\{s\}\ra C\times\{s\}$.

\

\subsection{Families of vector bundles in Case $U$}

\begin{prop}\label{families_U}
Let $\F_S$ be a rank $2$ vector bundle on $\PP^1\times C\times S$.
Assume that the bundle $\F_S$ has generic fiber type $(a,k-a)$
with $a>\frac{k}{2}$. Then there is a dense open $S^0\subseteq S$
with the following properties: \bi \item[i) ]For $s\in S^0$, the
restriction of the canonical extension of $\F_S$ to $\PP^1\times
C\times\{s\}$ is, up to multiplication by a scalar, the canonical
extension of $\F_s$ \item[ii) ]There is an integer $d$ such that
for $s\in S^0$ the canonical line subbundle $\F'_s$ of $\F_s$ has
type $(a,d)$. If $s_0\in S\setm S^0$, then $\F'_{s_0}$ has type
$(a_0,d_0)$ with $a_0\geq a$ and $d_0\geq d$. \ei
\end{prop}

\bp
Let $\eta$ be the generic point of $S$. The generic fiber type of $\F_{\eta}$
is $(a,k-a)$. Consider the canonical extension of $\F_{\eta}$ on
$\PP^1\times C\times \{\eta\}$:
\begin{equation}\label{can_seq_generic}
0\ra{\F'}_{\eta}\ra\F_{\eta}\ra\J_{\eta}\ra0
\end{equation}
where ${\F'}_{\eta}$ is a line bundle and $\J_{\eta}$ is a
torsion-free sheaf. Let ${\F'}_{\eta}\cong O(a)\bo\N_{\eta}$,
where $\N_{\eta}$ is a line bundle on $C\times\{\eta\}$.

We prove that  when we restrict the canonical extension of $\F_S$
to $\PP^1\times C\times \{\eta\}$ we get exactly the exact
sequence (\ref{can_seq_generic}). Let $\F_1\ra\F_S$ be an
injective morphism such that when restricted to $\PP^1\times
C\times \{\eta\}$ we get the injective morphism
$\F'_{\eta}\ra\F_{\eta}$ of (\ref{can_seq_generic}). Let
$\J_1=\coker(\F_1\ra\F_S)$. By restricting to $\PP^1\times C\times
\{\eta\}$, we get an exact sequence isomorphic to
(\ref{can_seq_generic}):
\begin{equation*}
0\ra(\F_1)_{\eta}\ra\F_{\eta}\ra(\J_1)_{\eta}\ra0
\end{equation*}

Note that since ${\F'}_{\eta}\cong(\F_1)_{\eta}$, we have
$\rk\F_1=1$. Let $\F'_S$ be the saturation of $\F_1$ in $\F_S$ and
let $\J_S$ be the torsion free sheaf which is the cokernel of the
morphism $\F_1\ra\F_S$. We prove that the exact sequence
\begin{equation}\label{can_seq_F_S}
0\ra\F'_S\ra\F_S\ra\J_S\ra0
\end{equation}
is the canonical extension of $\F_S$. Note that
$\rk\F'_S=\rk\F_1=1$. By Lemma \ref{johan}, it follows that
$\F'_S$ is a reflexive sheaf. As a reflexive sheaf of rank $1$ is
locally free, therefore we have that $\F'_S$ is a line bundle.
From the isomorphisms
\begin{equation*}
(\F'_S)_{\eta}\cong(\F_1)_{\eta}\cong\F'_{\eta}\cong O(a)\bo\N_{\eta}
\end{equation*}
it follows that $\F'_S\cong O(a)\bo\N_S$, for some line bundle $\N_S$ on
$C\times S$. By Lemma \ref{this_is_can_U}, the sequence
(\ref{can_seq_F_S}) is the canonical extension of $\F_S$, up to
multiplication by a scalar.


Let $s\in S$ and consider the canonical extension of $\F_s$:
\begin{equation*}
0\ra\F'_s\ra\F_s\ra\J_s\ra0
\end{equation*}
where $\F'_s\cong O(a_s)\bo\Lb_s$, for some integer $a_s$, and $\Lb_s$
is the line bundle ${p^s_2}_*\F_s(-a)$.

By restricting the exact sequence (\ref{can_seq_F_S}) to
$\PP^1\times C\times\{s\}$ we get an exact sequence:
\begin{equation}\label{can_seq_F_s}
0\ra(\F'_S)_s\ra\F_s\ra(\J_S)_s\ra0
\end{equation}

The map  $(\F'_S)_s\ra\F_s$ factors through the
canonical morphism $(\F'_S)_s\ra\F'_s$, which therefore must be injective.
Since $\F'_S\cong O(a)\bo\N_S$ and $\F'_s\cong O(a_s)\bo\Lb_s$, it follows
that  $a_s\geq a>\frac{k}{2}$ and there is an injective morphism
$(\N_S)_s\ra\Lb_s$.

\ps

Let $d$ be the integer such that $\deg(\N_S)_s=d$ for all $s\in
S$. It is clear now that for any $s\in S$, we have $\deg\Lb_s\geq
d$. We claim that there is a dense open $S^0\subset S$ such that
for $s\in S^0$ the sheaf $(\J_S)_s$ on $\PP^1\times C\times\{s\}$
is torsion free. Assuming the claim, we have by Lemma
\ref{this_is_can_U} that, for $s\in S^0$, the sequence
(\ref{can_seq_F_s}) is, up to multiplication by a scalar, the
canonical extension of $\F_s$. It follows that, for $s\in S^0$,
the map $(\F'_S)_s\ra\F_s$ is an isomorphism, therefore $a_s=a$
and $(\N_S)_s\cong\Lb_s$. In particular, for $s\in S^0$, we have
$\deg\Lb_s=d$.

\ps

We prove now the claim. Since $\J_S$ is torsion free of rank $1$,
there is a closed subscheme $\D$ of $\PP^1\times C\times S$, of
codimension at least $2$, and a line bundle $\M_S$, such that
$\J\cong\I_{\D}\otimes\M_S$,
where $\I_{\D}$ is the ideal sheaf of $\D$. 

If $s\in S$ is a point (not necessarily closed), let $\D_s$ be the
scheme theoretic intersection $\D\cap(\PP^1\times C\times\{s\})$ and
let $\T_s$ be the sheaf $\TTor_1(O_{\D},O_{\PP^1\times C\times\{s\}})$.
Note that $\T_s$ is supported on $\D_s$. There is an exact
sequence on $\PP^1\times C\times\{s\}$:
\begin{equation*}
0\ra\T_s\ra(\I_{\D})_s\ra\I_{\D_s}\ra0
\end{equation*}
where $\I_{\D_s}$ is the ideal sheaf of the closed subscheme
$\D_s$ in $\PP^1\times C\times\{s\}$. It follows that $\T_s$ is
the torsion subsheaf of $(\I_{\D})_s$. Hence, for any point $s\in
S$, we have $\T_s=0$ if and only if $(\I_{\D})_s$ (and hence
$(\J_S)_s$) is torsion-free.

Note that $\T_s=0$ if and only if $\D$ is flat over $S$ at the
points in $\D_s$. Let $S^0$ be the open set of those $s\in S$ for which
$\D$ is flat over $S$ at the points in $\D_s$.  By our construction,
$(\J_S)_{\eta}$ is torsion-free, therefore  $\T_{\eta}=0$ and $\D$
is flat over $S$ at the points in $\D_{\eta}$. It follows that
$S^0$ is dense in $S$. \ep

Recall that a sheaf $\F$ is called \emph{reflexive} if the canonical
map $\F\ra\F^{**}$ is an isomorphism and that the sheaf $\F$ is torsion-free
if and only if $\F\ra\F^{**}$ is injective.

\begin{lemma}\label{johan}
Let $0\ra\F'\ra\F\ra\F''\ra0$ be an exact sequence. If $\F$ is reflexive
and $\F''$ is torsion free, then $\F'$ is reflexive.
\end{lemma}

\bp Since $\F$ is reflexive, in particular, it is torsion-free. It
follows that $\F'$ is torsion-free; hence the canonical map
$\al:\F'\ra{\F'}^{**}$ is injective. One may find an inverse to
$\al$, by using that the composition
${\F'}^{**}\ra{\F}^{**}\ra{\F''}^{**}$ is zero and the fact that
$\F\ra{\F}^{**}$ is an isomorphism. \ep

\begin{lemma}\label{this_is_can_U}
Let $T$ be an irreducible variety and $\F$ a rank $2$
vector bundle on $\PP^1\times T$ of fiber degree $k$.  Assume there is
an integer $a>\frac{k}{2}$ and a short exact sequence on $\PP^1\times T$:
\begin{equation}\label{U_can_looking}
0\ra O(a)\bo\Lb\ra\F\ra\J\ra0.
\end{equation}
where $\Lb$ is a line bundle on $T$ and $\J$ is a torsion free
sheaf on $\PP^1\times T$. 
Then the extension (\ref{U_can_looking}) is, up to scalar
multiplication, the canonical extension of $\F$.
\end{lemma}
\bp Assume that $\F$ has generic fiber type $(b,k-b)$, with
$b\geq\frac{k}{2}$. If $\eta$ is the generic point of $T$,  by restricting
(\ref{U_can_looking}) to $\PP^1\times\{\eta\}$, we get that
$\J_{\eta}\cong O(k-a)$ and $b=a$. 

Tensor (\ref{U_can_looking}) by $O(-a)$ and take the push-forward to $T$:
\begin{equation*}
0\ra\Lb\ra{p_2}_*\F(-a)\ra{p_2}_*\J(-a)\ra0
\end{equation*}

Note that ${p_2}_*\J(-a)_{\eta}\cong\H^0(\PP^1_{\eta},O(k-2a)=0$.
It follows that  ${p_2}_*\J(-a)$ is a torsion sheaf. But as $\J$
is torsion-free, ${p_2}_*\J(-a)$ is torsion-free, so
${p_2}_*\J(-a)=0$. It follows that $\Lb\cong{p_2}_*\F(-a)$. As
${p_2}_*\F(-a)_{\eta}\cong\H^0(\PP^1_{\eta},O\oplus
O(k-2a))$ and since ${p_2}_*\F(-a)\cong\Lb$ has rank $1$, we have
$b=a$.

We proved that the canonical subbundle $\F'$ of $\F$ is isomorphic to
$O(a)\bo\Lb$. Since there are no non-zero maps $\F'\ra\J$,
there is an induced commutative diagram between (\ref{U_can_looking}) and
the canonical sequence of $\F$.
Since a non-zero endomorphism of $\F'$ is given by scalar multiplication,
the two sequences are scalar multiples of each other.
\ep


\subsection{Families of vector bundles in Case $E$}\label{fam_E}

\

\ps

In this section we will work with families of vector bundles on
$\PP^1\times C$ with equal generic splitting $(a,a)$. Since we can
tensor $\F$ with $O(-a)$, we may assume that $a=0$.

\begin{prop}\label{families_E}
Let $\F_S$ be a rank $2$ vector bundle on $\PP^1\times C\times S$
such that, for any $s\in S$, the generic fiber type of the bundle
$\F_s$ is $(0,0)$. Assume that for any $s$ in $S$ and for any
$c\in C$ the bundle $\F_{s,c}$ splits either as $O\oplus O$ or
$O(1)\oplus O(-1)$. Then for any $s\in S$ we have: \bi \item[i)
]The canonical extension of $\F_s$ has the form:
\begin{equation}\label{can_ext_F_s} 
0\ra O\bo\E_s\ra\F_s\ra O(-1)\bo O_{D_s}\ra0
\end{equation}
where $\E_s$ is a rank $2$ vector bundle on $C$ and $D_s$ is a
$0$-cycle on $C$ \item[ii) ]The extension (\ref{can_ext_F_s}) and
the restriction of the canonical extension of $\F_S$ to
$\PP^1\times C\times\{s\}$ are in the same orbit for the action of
$\Aut(O_{D_s})$ on $\Ext^1(O(-1)\bo O_{D_s},O\bo\E_s)$. \ei In
particular, $\deg(\E_s)$ is constant for any $s\in S$.
\end{prop}

\bp
Consider the canonical extension of $\F_S$:
\begin{equation}\label{can_ext_F_S}
0\ra\F'_S\ra\F_S\ra\J_S\ra0
\end{equation}

Recall that $\F'_S$ is the rank $2$ torsion free sheaf
$({p_{2,3}}_*\F_S)$, where $p_{2,3}$ is the canonical projection
$\PP^1\times C\times S\ra C\times S$. For any $s\in S$, there is a
canonical morphism $\phi_s:(\F'_S)_s\ra\F'_s$. We prove that
$\phi_s$ is an isomorphism. For any $c\in C$ and $s\in S$ there
are canonical morphisms
\begin{equation*}
\psi_{c,s}:({p_{2,3}}_*\F)_{c,s}\ra\H^0(\PP^1,\F_{c,s}),\quad
\chi_{c,s}:({p^s_2}_*\F_s)_c\ra\H^0(\PP^1,\F_{c,s})
\end{equation*}
such that $\chi_{c,s}\circ\phi_{c,s}=\psi_{c,s}$, where
$\phi_{c,s}$ is the restriction of the morphism $\phi_s$ to
$\PP^1\times\{c\}\times\{s\}$ (see \cite{M}). Since for
any $s\in S$ and $c\in C$ we have $\H^1(\PP^1,\F_{c,s})=0$, it
follows that $\psi_{c,s}$ and $\chi_{c,s}$ are isomorphisms.
Hence, $\phi_{c,s}$ is an isomorphism for any $c\in C$ and
therefore, we have that $\phi_s$ is an isomorphism.

Let $s\in S$ and consider the canonical extension of $\F_s$:
\begin{equation*}
0\ra\F'_s\ra\F_s\ra\J_s\ra0
\end{equation*}

Since the bundle $\F_s$ has generic type $(0,0)$ we have
$\F'_s\cong O\bo\E_s$, where $\E_s={p^s_2}_*\F_s$ is a rank
$2$ bundle. 
Since for any $s$ in $S$ and for any $c\in C$ the
bundle $\F_{s,c}$ splits either as $O\oplus O$ or $O(1)\oplus O(-1)$, it
follows that $\J_s\cong O(-1)\bo O_{D_s}$, where $D_s$ is some
$0$-cycle on $C$.

The restriction map $(\F'_S)_s\ra\F_s$ factors through the
canonical morphism $\phi_s:(\F'_S)_s\ra\F'_s$, which we proved is
an isomorphism. It follows that the restriction map
$(\F'_S)_s\ra\F_s$ is injective and by restricting
(\ref{can_ext_F_S}) to $\PP^1\times C\times\{s\}$ we get an exact
sequence:
\begin{equation}
0\ra(\F'_S)_s\ra\F_s\ra(\J_S)_s\ra0
\end{equation}

Moreover, there is an isomorphism $(\J_S)_s\cong\J_s$.
Since $\Aut(\J_s)\cong\Aut(O_{D_s})$, the proposition follows.
\ep

\begin{cor}\label{stable_in_fam}
Under the conditions and with the notations in Proposition
\ref{families_E}, there is an open $S^0\subseteq S$ (possibly
empty) such that for any $s\in S^0$, the bundle $\E_s$ is a
\emph{stable} bundle.
\end{cor}

\bp
Let $\F'_S$ be the canonical subbundle of $\F_S$. Then we have that
\begin{equation*}
\E_s={p^s_2}_*\F_s\cong({p_{2,3}}_*\F)_s
\end{equation*}

Since being stable is an open condition, the locus of those $s\in
S$ such that $({p_{2,3}}_*\F)_s$ is stable is open in $S$. It
follows that there is an open $S^0\subseteq S$ (possibly empty)
such that for any $s\in S^0$, the bundle $\E_s$ is stable. \ep

\

\textbf{Comment. }
Note that for any rank $2$ vector bundle $\F_S$ on $\PP^1\times C\times S$
with generic fiber type $(0,0)$, we have by upper-semicontinuity that
there is a dense open $S^0\subset S$ such that for any $s\in S^0$,
the bundle $\F_s$ has generic fiber type $(0,0)$. Therefore,
the first assumption
in Proposition \ref{families_E} is not hard to satisfy. However,
the assumption about the splitting of $\F_{s,c}$ is a non-trivial one
(see Remark \ref{complete}).

\

\textbf{Complete families of vector bundles}

\

If $T$ and $X$ are \emph{smooth} irreducible
varieties, a vector bundle $\F$ on $X\times T$ is called a \emph{complete}
family of vector bundles on $X$ if the Kodaira-Spencer infinitesimal
deformation map
\begin{equation*}
\om:T_t T\ra\Ext^1_X(\F_t, \F_t)
\end{equation*}
is surjective for any $t\in T$ \cite{P} (II.15.1) .

If $X=\PP^1$ and $\F$  is a complete family of rank $2$ vector bundles
on $\PP^1$, such that for general $t\in T$ we have that $\F_t\cong O\oplus O$,
then the points $t\in T$ for which $\F_t$ splits as
$O(\al)\oplus O(\be)$, with  $|\al-\be|\geq3$ form a closed set of
codimension at least $2$ (\cite{P}, II.15.4.3). By applying this to 
$T=C\times S$, where $S$ is a smooth variety,
we have the following consequence. 

\begin{rmk}\label{complete}
Let $\F_S$ be a rank $2$ vector bundle on $\PP^1\times C\times S$,
with generic fiber type $(0,0)$. If $S$ is smooth and $\F$ is a
complete family of vector bundles on $\PP^1$, there is a dense
open $S^0\subset S$ such that for $s\in S^0$ and $c\in C$, the
bundle $\F_{s,c}$ splits either as $O\oplus O$ or $O(1)\oplus
O(-1)$.
\end{rmk}

Note that if $\F_S$ is a vector bundle on $\PP^1\times C\times S$,
then  $\F_S$ can be regarded as a family of vector bundles on $\PP^1\times C$,
as well as a family of vector bundles on $\PP^1$.
Consider the Kodaira-Spencer map given by the
family $\F_S$ of vector bundles on $\PP^1\times C$:
\begin{equation}\label{KS_global}
\om:T_s S\ra\Ext^1_{\PP^1\times C}(\F_{s}, \F_{s})
\end{equation}

If $\F_S$ is such that for any $s\in S$ and $p\in\PP^1$ we have
$\det(\F_{s,p})\cong\xi$, the Kodaira-Spencer map (\ref{KS_global})
factors through a map:
\begin{equation}\label{KS_global+}
\om^+:T_s S\ra\Ext^1_{\PP^1\times C}(\F_{s}, \F_{s})^+
\end{equation}
where $\Ext^1_{\PP^1\times C}(\F_{s}, \F_{s})^+$ is the subspace
in $\Ext^1_{\PP^1\times C}(\F_{s}, \F_{s})$ corresponding to traceless
endomorphisms. Via the canonical isomorphism 
\begin{equation*}
\Ext^1_{\PP^1\times C}(\F_{s}, \F_{s})\cong
\H^1(\PP^1\times C,\F_s\otimes\F_s^*)
\end{equation*}
the space $\Ext^1_{\PP^1\times C}(\F_{s}, \F_{s})^+$ maps to
$\H^1(\PP^1\times C,\ad(\F_s))$, where by $\ad(\F_s)$ we denote
the subbundle of $\F_s\otimes\F_s^*$ of endomorphisms of trace zero.

\ps

If we fix a closed point $c\in C$, the family 
$\F_c={\F_S}_{|\PP^1\times\{c\}\times S}$, of
vector bundles on $\PP^1$, induces a Kodaira-Spencer map:
\begin{equation}\label{KS_local}
\om_c:T_s S\ra\Ext^1_{\PP^1}(\F_{s,c}, \F_{s,c})
\end{equation}

\ps

It is easy to check, using the definition of the Kodaira-Spencer
map (see \cite{P}, II.15.1) that $r\circ\om=\om_c$, where
$r$ is the restriction map
\begin{equation}\label{res_c}
r:\Ext^1_{\PP^1\times C}(\F_s,\F_s)\ra\Ext^1_{\PP^1}(\F_{s,c},\F_{s,c})
\end{equation}

The following lemma will help us to check when a rank $2$ vector
bundle $\F_S$ on $\PP^1\times C\times S$ is a complete family of
vector bundles on $\PP^1$.
\begin{lemma}\label{global_local_complete}
Let $\F_S$ be a rank $2$ vector bundle on $\PP^1\times C\times S$,
with generic fiber type $(0,0)$. Assume that for any $s\in S$ and 
$p\in\PP^1$ we have $\det(\F_{s,p})\cong\xi$. If the map $\om^+$ in
(\ref{KS_global+}) induced by $\F_S$ is surjective for any $s\in S$, 
then $\F_S$ is a complete family of vector bundles on $\PP^1$.
\end{lemma}

\bp
Let $s\in S$, $c\in C$. 
By Lemma \ref{restr_is_surj}, the restriction $r^+$
of the map $r$ of (\ref{res_c}) to $\Ext^1_{\PP^1\times C}(\F_s,\F_s)^+$
is surjective. It follows that the composition
$\om_c=\om\circ r$ is surjective. Hence, the family $\F_c$, of vector
bundles on $\PP^1$, is complete, for any $c\in C$.

Consider the Kodaira-Spencer map induced by $\F_S$ as a family of vector
bundles on $\PP^1$:
\begin{equation*}
\ups:T_s S\oplus T_c C\ra\Ext^1_{\PP^1}(\F_{s,c}, \F_{s,c})
\end{equation*}

It is a straightforward computation to check that,
if $i$ is the inclusion $T_s S\ra T_s S\oplus T_{c} C$ given by
$u\mapsto(u,0)$, then the composition $\ups\circ i$ is the Kodaira-Spencer
map $\om_c$, given by the family $\F_c$. It follows that the family $\F_S$
of vector bundles on $\PP^1$ is complete.
\ep

\begin{lemma}\label{restr_is_surj}
Let $\G$ be a rank $2$ vector bundle on $\PP^1\times C$ with generic fiber
type $(0,0)$. For any $c\in C$, the restriction map 
$\Ext^1_{\PP^1\times C}(\G,\G)\ra\Ext^1_{\PP^1}(\G_c, \G_c)$
induces a surjective map
\begin{equation*}
\Ext^1_{\PP^1\times C}(\G,\G)^+\ra\Ext^1_{\PP^1}(\G_c, \G_c)
\end{equation*}
\end{lemma}

\bp
It is enough to prove that
$\H^2(\PP^1\times C,\ad(\G)\otimes\I_c)=0$, where $\I_c$ is the
ideal sheaf of $\PP^1\times\{c\}$ in $\PP^1\times C$.
Using a Leray spectral sequence, we have
\begin{equation*}
\H^2(\PP^1\times C,\ad(\G)\otimes\I_c)\cong
\H^1(C,R^1{p_2}_*(\ad(\G)\otimes\I_c))
\end{equation*}

Let $\ze$ be the generic point of $C$. Then $\G_{\ze}\cong O\oplus O$
and $R^1{p_2}_*(\ad(\G)\otimes\I_c)_{\ze}=0$. It follows that
the sheaf  $R^1{p_2}_*(\ad(\G)\otimes\I_c)$ is supported
on a finite set of points of $C$. It follows:
\begin{equation*}
\H^1(C,R^1{p_2}_*(\ad(\G)\otimes\I_c))=0
\end{equation*}
\ep


\subsection{Proof of Theorem \ref{M^0}}\label{pf_M^0}

\

\ps

Let $\M\subset\Hom_k(\PP^1,M)$ be an irreducible component and let
$\HH_{\M}$ be the restriction of the universal bundle $\HH$ to
$\PP^1\times C\times\M$. For $f\in\M$, we have that the bundle
$\HH_f$ on $\PP^1\times C$ is isomorphic to $(f\times{id})^*\U$. By
Lemma \ref{corresp}, the fiber degree of the bundle $\HH_f$ is $k$
and the Chern classes of $\HH_f$ are $c_1$ and $c_2$ as in
(\ref{c}). Let $a\geq\frac{k}{2}$ be the integer
such that $\HH_{\M}$ has generic splitting type $(a,k-a)$.

\textbf{Case $a>\frac{k}{2}$. }

By Proposition \ref{families_U}, it follows that there is an
integer $e$ and a dense open $\M^0\subset\M$ such that for
$f\in\M^0$, the canonical subbundle of $\HH_f$ has type $(a,-e)$.
Since the bundle $\HH_f$ induces a morphism $f:\PP^1\ra M$, the
bundle $\HH_f$ corresponds to a point of the good of the moduli
space $\Bb(a,e)$. Hence, the members of the family $\HH_{\M}$
correspond to points in $\Bb^0(a,e)$.

\ps

\textbf{Case $a=\frac{k}{2}$. }

Assume $k=2a$ is an even integer. Let
$\M\subset\Hom_k(\PP^1,M)$ be an irreducible component with
generic fiber type $(a,a)$ and let $\M^0\subset\M$ be a
dense open such that for $f\in\M^0$ the vector
bundle $\HH_f$ has generic fiber type $(a,a)$. By \cite{C} (Lemma
3.4), $f$ is a smooth point of $\Hom_k(\PP^1,M)$. We have that
$\M$ is an irreducible component of the expected dimension
$(2k+3g-3)$, which is smooth along $\M^0$. 

\ps

For $f\in\M^0$, using that $\Hom_k(\PP^1,M)$ is smooth along $\M^0$, 
we have that the family $\HH_{\M^0}$ of vector 
bundles on $\PP^1\times C$ induces a Kodaira-Spencer map
\begin{equation}\label{KS}
\om^+:T_f\Hom_k(\PP^1,M)\ra\Ext^1_{\PP^1\times C}(\HH_f, \HH_f)^+\cong
\H^1(\PP^1\times C, \ad(\HH_f))
\end{equation}

We prove this map is an isomorphism. First, recall that
there is a canonical isomorphism
\begin{equation}\label{canonic_1}
T_f\Hom_k(\PP^1,M)\ra\H^0(\PP^1,f^*T_M)
\end{equation}
(See \cite{K}, II.1.7). On the other hand, 
we have $f^*T_M\cong\R^1{p_1}_*\ad(\HH_f)$ \cite{N2}.
One may check directly, using deformation theory,
that if we compose (\ref{KS}) with the canonical map
\begin{equation}\label{canonic_2}
\H^1(\PP^1\times C, \ad(\HH_f))\ra\H^0(\PP^1, \R^1{p_1}_*\ad(\HH_f))
\end{equation}
we obtain the map (\ref{canonic_1}). Since ${p_1}_*\ad(\HH_f)=0$,
the canonical map (\ref{canonic_2}) is an isomorphism.
It follows that the  Kodaira-Spencer map (\ref{KS}) is an isomorphism.
 
\ps

By Lemma \ref{global_local_complete}, the family $\HH_{\M^0}$ of vector
bundles on $\PP^1$ is complete. Using Remark \ref{complete}, it
follows that, eventually shrinking $\M^0$, if $f\in\M^0$, for any
$c\in C$, the bundle $\HH_{f,c}$ splits either as $O(a)\oplus
O(a)$ or as $O(a+1)\oplus O(a-1)$. By Proposition
\ref{families_E}, it follows that, for any $f\in\M^0$, the
canonical extension of the bundle $\HH_f$ has the form
\begin{equation*}
0\ra O(a)\bo\E\ra\HH_f\ra O(a-1)\bo O_{D}\ra0
\end{equation*}
for some rank $2$ vector bundle $\E$  on $C$ and $D$ a $0$-cycle on
$C$. By a computation with Chern classes, it follows that
$\deg(D)=a$ and $\deg(\E)=(1-a)$.

We claim that, for general $f\in\M^0$, the bundle $\E$ is stable.
This proves that (eventually shrinking $\M^0$) there is a dense
open $\M^0\subset\M$, such that for any $f\in\M^0$, the bundle
$\HH_f$ has type $(\da)$ (see (\ref{type_da})), and so the  
members of the family $\HH_{\M^0}$ correspond to points in $\Bbe^0$.

We prove now the claim. By Corollary \ref{stable_in_fam}, it follows that
there is an open in $\M^0$, possibly empty, such that the bundle
$\E$ is stable. We prove that this open is dense. Let
$f\in\Hom_k(\PP^1,M)$ be such that the bundle $\E$ is not stable.
Then there is a line subbundle $\Lb'$ of $\E$, with
$\deg(\Lb')\geq\frac{1-a}{2}$. Let $\F'$ be the saturation of
$O(a)\bo\Lb'$ in $\HH_f$ and let $\J$ be the torsion free
quotient:
\begin{equation}\label{seq_H_f}
0\ra\F'\ra\HH_f\ra\J\ra0
\end{equation}

Since $\HH_f$ is locally free, by Lemma \ref{johan}, the sheaf
$\F'$ is reflexive. Since a reflexive sheaf of rank $1$ is locally
free, we have that $\F'$ is a line bundle. We have $\F'\cong
O(a)\boxtimes\Lb$, for some line bundle $\Lb$ on $C$. If we let
$d=-\deg\Lb$, note that $d\leq\frac{a-1}{2}$. Since $\HH_{f}$ is
stable on the fibers over $\PP^1$, we have $d\geq0$.

As $\J$ is torsion free, we have $\J\cong\I_Z\otimes\M$, for some
$Z\subset\PP^1\times C$ a $0$-cycle and $\M$ a line bundle on $\PP^1\times C$.
From a computation with Chern classes, we have $\M\cong
O(a)\bo(\Lb^{-1}\otimes\xi)$.

By Remark \ref{extra_B}, there exists a variety $\Bb(a,d)$ which
parameterizes extensions (\ref{seq_H_f}). Moreover, since there is
a universal bundle on $\PP^1\times C\times\Bb$, by the same
argument as in Remark \ref{T'}, there is a morphism
$\mu_d:\Bb^0(a,d)\ra\Hom_k(\PP^1,M)$. It follows that $f$ is in
the image of the morphism $\mu_d$, for some $0\leq
d\leq\frac{a-1}{2}$.

Note that the dimension of $\Bb(a,d)$ may be computed from
(\ref{dim_B(a,e)}) to be:
\begin{equation*}
\dim\Bb(a,d)=3a+2g+2d-1
\end{equation*}

For any $d$ such that $0\leq d\leq\frac{a-1}{2}$, it follows that
$\dim\Bb(a,d)<(4a+3g-3)=(2k+3g-3)$. It follows that a general
$f\in\M$ is not in the image of $\mu_d$ for any $d$.
Hence, for general $f\in\M$, the bundle $\E$ in the canonical
sequence of $\HH_f$ is stable.


\section{Irreducible components coming from $\M(a,e)$}

Fix $k\geq1$ and consider the range ($\star$) for the pairs of
integers $(a,e)$ for which the good locus $\Bb(a,e)$ is non-empty (see Figure 1):
\begin{equation*}\tag{$\star$}\label{star}
\{ (a,e)\quad|\quad k\geq a>k/2,\quad\frac{k-a}{2a-k}\geq
e>0\}\cup\{(k,0)\}
\end{equation*}
Let $A$ be the branch of the hyperbola in the plane with axis $a$
and $e$ given by:
\begin{equation*}\tag{$A$}
e=\frac{k-a}{2a-k}=-\frac{1}{2}+\frac{k}{2(2a-k)}, \quad
a>\frac{k}{2}
\end{equation*}

\begin{figure}[h]
\centerline{ \psfig{figure=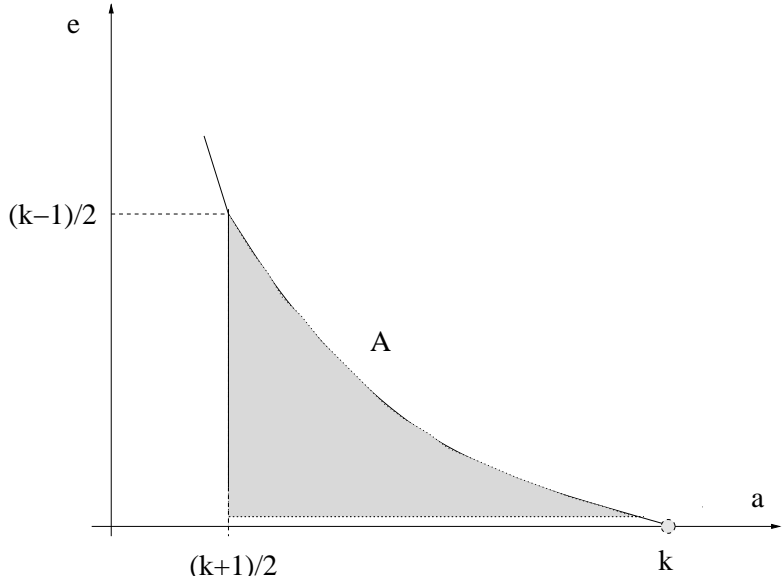,height=3in} }
\smallskip
\centerline{
\parbox{4.5in}{\caption[Range ($(\star)$)]{{\small This shows the
range $(\star)$ of the pairs of integers $(a,e)$. The curve $A$ has
equation $e=\frac{k-a}{2a-k}$.}}} }
\end{figure}

For $(a,e)$ in the range (\ref{star}), the subvarieties $\M(a,e)$
of $\Hom_k(\PP^1, M)$ have dimension
\begin{equation}\label{dim}
\dim\M(a,e)=(2a-k+2)g+(3k-3a-1)-e(2a-k-2)
\end{equation}

\begin{obs}\label{dim_grow}
For integers $a$ and $e$ in the range (\ref{star}), the dimensions of
the subvarieties  $\M(a,e)$ grow as follows (see Figures 2 and 3):
\bi
\item[i. ]If $a=\frac{k+1}{2}$ and $e'>e$, then
$\dim\M(a,e')>\dim\M(a,e)$
\item[ii. ]If $a=\frac{k}{2}+1$, then for any $e$,
$\dim\M(a,e)=4g+\frac{3k}{2}-4$
\item[iii. ]If $k\geq a>\frac{k}{2}+1$ and  $e'>e$, then
$\dim\M(a,e')<\dim\M(a,e)$
\item[iv. ]If $a>a'$ and $e\geq (g-1)$, then
$\dim\M(a,e)<\dim\M(a',e)$
\item[v. ]If $a>a'$ and $e\leq (g-2)$, then
$\dim\M(a,e)>\dim\M(a',e)$
\ei
\end{obs}

\begin{figure}[h]
\centerline{
\psfig{figure=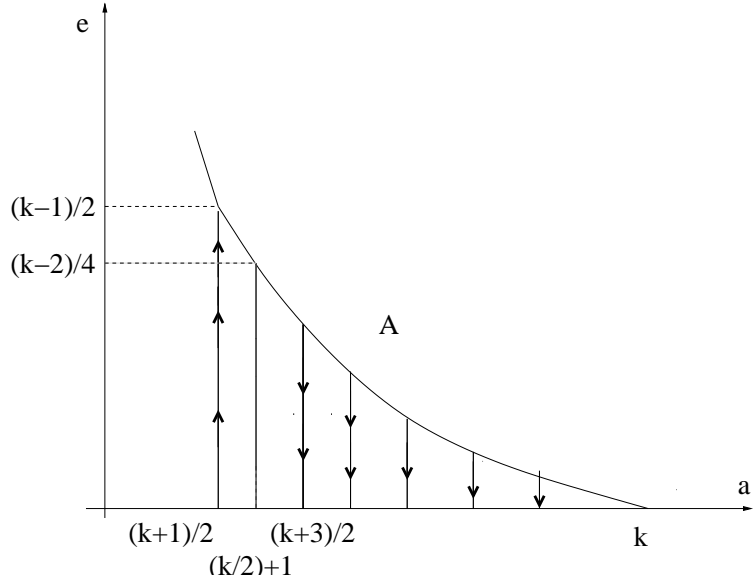,height=3in}
}
\smallskip
\centerline{
\parbox{4.5in}{\caption[]{{\small The arrows show the
direction in which dimensions grow when $a$ is constant. If $k$ is
even and $a=\frac{k}{2}+1$, then the dimension is constant for any
$e$ in the given range.}}}
}
\end{figure}

\begin{figure}[h]
\centerline{
\psfig{figure=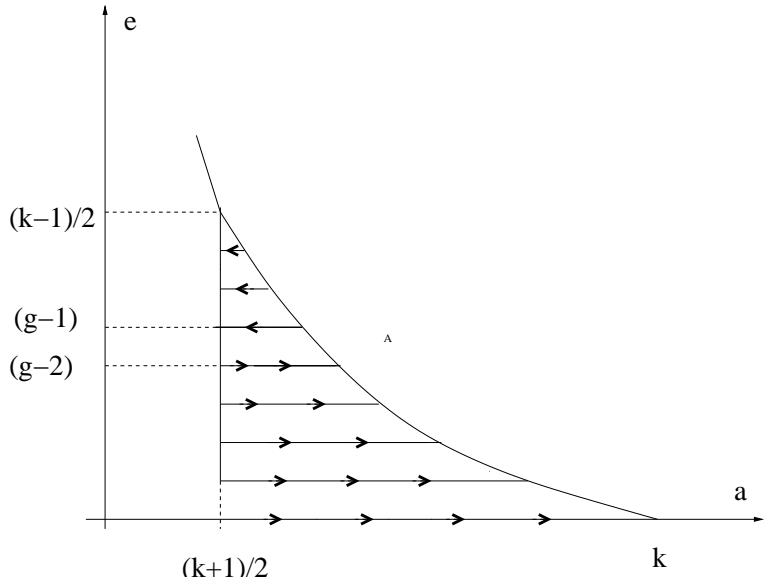,height=3in}
}
\smallskip
\centerline{
\parbox{4.5in}{\caption[]{{\small The arrows show the
direction in which dimensions grow when $e$ is constant. In this
picture we assumed $\frac{k-1}{2}\geq(g-1)$. }}}
}
\end{figure}

\bp
This is clear since
\begin{gather*}
\dim\M(a,e')-\dim\M(a,e)=(2a-k-2)(e-e')\\
\dim\M(a',e)-\dim\M(a,e)=(2g-2e-3)(a'-a).
\end{gather*}
\ep

Note that if $\Mae$ is an irreducible component, then
$\dim\M(a,e)\geq(2k+3g-3)$. The locus of pairs $(a,e)$ in the
range (\ref{star}) where $\dim\M(a,e)=(2k+3g-3)$ is given by: \bi
\item[i) ] The point $(a,e)=(\frac{k+1}{2},\frac{k-1}{2})$ 
\item[ii) ]On the line $a=\frac{k}{2}+1$, if $k=2g-2$ 
\item[iii) ]The part of the hyperbola $B$, given by equation:
\begin{equation*}\tag{$B$}\label{B}
e=\frac{2g-3}{2}+\frac{2g-2-k}{2(2a-k-2)}
\end{equation*}
that is inside the range (\ref{star}), when $a\geq\frac{k+3}{2}$.
\ei

Note that if $k=2g-2$, then $B$ is the line $e=\frac{2g-3}{2}$.

\begin{defn}
Define $R$ to be the part of the region (\ref{star}) where the
dimension of $\M(a,e)$ is at least the expected one, that is, one
of the following: \bi \item[i) ]The point
$(\frac{k+1}{2},\frac{k-1}{2})$, if $k$ odd \item[ii) ]The part of
the line $a=\frac{k}{2}+1$ inside the region (\ref{star}), when
$k$ is even and $k\leq(2g-2)$ \item[iii) ]The part of the region
between the graphs of $A$ and $B$, where the dimension is at least
the expected one, when $a\geq\frac{k+3}{2}$ (see Figure 4 and for 
special cases see the Figures in Section \ref{region_R}) 
\ei
\end{defn}

\begin{figure}[h]
\centerline{
\psfig{figure=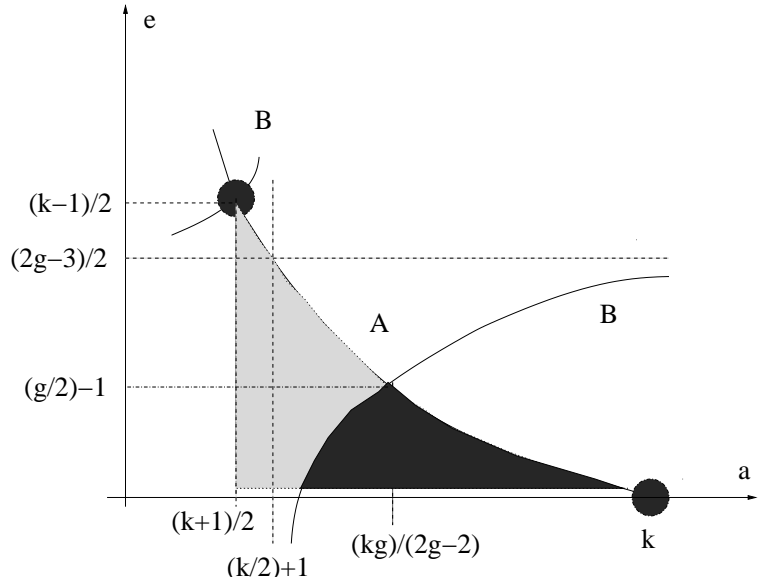,height=3in}
}
\smallskip
\centerline{
\parbox{4.5in}{\caption[]{{\small The dark shaded part is the region $R$
when $g\geq3$ and $k>2g-2$. The
hyperbola $B$ is the locus where the dimension is equal to the
expected dimension}}} }
\end{figure}

The hyperbola $B$ intersects the curve $A$ at the points:
\begin{equation*}
a=\frac{kg}{2(g-1)},\quad e=\frac{g}{2}-1\quad\hbox{and}\quad
a=\frac{k+1}{2},\quad e=\frac{k-1}{2}
\end{equation*}

The curve $B$ intersects the $a$-axis at the point
\begin{equation*}
a=\frac{(k+1)(g-1)-1}{2g-3}
\end{equation*}

The following observation is a consequence of Proposition
\ref{families_U}.
\begin{obs}\label{specialize}
For each $(a,e)$ in (\ref{star}), the subvarieties $\M(a,e)$
satisfy the property that if $\M(a,e)\subset\overline{\M(a',e')}$,
then $a'\leq a$ and $e'\geq e$
\end{obs}

We recall as observations the following statements proved in
Section 4.
\begin{obs}\label{odd_n_c}
If $k$ is odd, then $\overline{\M(\frac{k+1}{2},\frac{k-1}{2})}$ is the
nice component.
\end{obs}

\begin{obs}\label{all_Mae}
If $\M$ is an irreducible component which is not the
nice component, there are uniquely determined $a$ and $e$ in the
range (\ref{star}), such that $\M=\Mae$
\end{obs}

\begin{thm}\label{main}
If a pair of integers $(a,e)$ is in the range (\ref{star}), then the
subvariety $\Mae$ is an irreducible component of $\Hom_k(\PP^1,M)$
if and only if $\M(a,e)$ has dimension bigger or equal than
expected, or, equivalently, the pair $(a,e)$ is in the region $R$
(see Figure 4).
Moreover: \bi \item[i) ]If the pair $(a,e)$ is not in the region
$R$, then $\Mae$ is contained in the nice component \item[ii)
]Together with the nice component, these are all the irreducible
components \ei
\end{thm}

\bp The condition is clearly necessary. We prove that it is also
sufficient. For $(a,e)=(\frac{k+1}{2},\frac{k-1}{2})$ (if $k$ is
odd), by Observation \ref{odd_n_c}, we have that $\Mae$ is an
irreducible component.

Let $(a,e)$ be in the region $R$ and
$(a,e)\neq(\frac{k+1}{2},\frac{k-1}{2})$. Assume that $\Mae$ is
not an irreducible component. Then there is
$\M'\subset\Hom_k(\PP^1,M)$ irreducible component containing
$\M(a,e)$ and we have $\dim\M'>\dim\M(a,e)\geq(2k+3g-3)$. It
follows that $\M'$ is not the nice component. By Observation
\ref{all_Mae}, there are integers $a'$ and $e'$ such that
$\M'=\overline{\M(a',e')}$. By Observation \ref{specialize}, it
follows that $a\geq a'$ and $e\leq e'$.

Note that in the region $R$ we have $e\leq\frac{g}{2}-1\leq(g-2)$ 
(see Figure 4 and also Figures in Section \ref{region_R}).
By Observation \ref{dim_grow} i., we have $a\geq\frac{k}{2}+1$, as on the 
line $a=\frac{k+1}{2}$ the only value of $e$ for which the dimension of 
$\M(a,e)$ is at least the expected one, is for $e=\frac{k-1}{2}$.
By Observation \ref{dim_grow}, we have
\begin{equation*}
\dim\M(a,e)\geq\dim\M(a',e)\geq\dim\M(a',e')
\end{equation*}
But since $\M'\neq\Mae$, we have $\dim\M(a',e')>\dim\M(a,e)$, which
gives a contradiction. It follows that  $\M(a,e)$ is an irreducible
component.

We prove now i). Let $(a,e)$ be in the range (\ref{star}), but not in the 
region $R$. The dimension of  $\M(a,e)$ is smaller than expected. Let
$\M'\subset\Hom_k(\PP^1,M)$ be an irreducible component, such that
$\M(a,e)\subset\M'$. Assume $\M'$ is not the nice component. By
Observation \ref{all_Mae}, there are integers $a'$ and $e'$ such
that $\M'=\overline{\M(a',e')}$. Since $\M'$ has dimension at
least the expected one, the pair $(a',e')$ is in the region $R$.
But by Observation \ref{specialize} it follows that $a\geq a'$ and
$e\leq e'$.

The only pair of integers in the range $R$ with this property is
$(\frac{k+1}{2},\frac{k-1}{2})$ which corresponds to the nice
component if $k$ is odd, which contradicts the assumption that  $\M'$
is not the nice component.
It follows that $\M'$ must be the nice component.

Part ii) is immediate from Observation \ref{all_Mae}.
\ep


\section{Description of the varieties $\M(a,e)$ and $\M_E$}

\subsection{MRC fibrations of the varieties $\M(a,e)$}

\begin{thm}\label{MRC_M(a,e)}
Let $(a,e)$ be a pair of integers in the range (\ref{star}) and
let $\de=(k-a)-e(2a-k)$. We distinguish the following cases: \bi
\item[i) ] If $\de=0$, the MRC fibration of the variety $\M(a,e)$
is given by a surjective map $$\M(a,e)\ra\Pic^{-e}(C)$$ \item[ii)
]If $\de>0$, the MRC fibration of the variety $\M(a,e)$ is given
by a map $$\M(a,e)\ra\Pic^{-e}(C)\times\Pic^{\de}(C) $$ which is
dominant if and only if $\de\geq g$. \ei
\end{thm}

\bp By Corollary \ref{M_isom_B^0}, there is an isomorphism
$\M(a,e)\cong\Bb^0(a,e)$. Recall that $\Bb^0(a,e)$ is a dense open
in a projective bundle over
$\Pic^{-e}(C)\times\Hilb^{\de}(\PP^1\times C)$; let $\pi$ be the
induced map:
\begin{equation*}
\pi:\M(a,e)\ra\Pic^{-e}(C)\times\Hilb^{\de}(\PP^1\times C)
\end{equation*}

The map $\pi$ associates to $f\in\M(a,e)$ the pair $(\Lb,Z)$ from
the canonical extension of the bundle $\HH_f$ corresponding to
$f$. Note that the fibration given by $\pi$ is locally trivial
with rational fibers. If $\de=0$, since $\Pic^{-e}(C)$ is an
abelian variety, the morphism $\pi$ gives the MRC fibration of
$\M(a,e)$.

If $\de>0$, there is a birational morphism
$v:\Hilb^{\de}(\PP^1\times C)\ra\Sym^{\de}(\PP^1\times C)$. Note
that
\begin{equation*}
\Sym^{\de}(\PP^1\times C)\cong
\Sym^{\de}(\PP^1)\times\Sym^{\de}(C)\cong
\PP^{\de}\times\Sym^{\de}(C)
\end{equation*}

Let  $u:\Sym^{\de}(C)\ra\Pic^{\de}(C)$ be the canonical morphism.
If $p_2:\PP^{\de}\times\Sym^{\de}(C)\ra\Sym^{\de}(C)$ is the second 
projection, consider the following composition, which we denote by $\si$:
\begin{equation*}
\begin{CD}
\Hilb^{\de}(\PP^1\times C)@>v>>\PP^{\de}\times\Sym^{\de}(C)@>{p_2}>>
\Sym^{\de}(C)@>u>>\Pic^{\de}(C)
\end{CD}
\end{equation*}

Since  $v$ is birational and the fibers of $u$ are rational,
it follows that the general fiber of $\si$ is rational. Denote by $\rho$ 
the following composition:
\begin{equation*}
\begin{CD}
\M(a,e)@>{\pi}>>\Pic^{-e}(C)\times\Hilb^{\de}(\PP^1\times
C)@>{\id\times\si}>> \Pic^{-e}(C)\times\Pic^{\de}(C)
\end{CD}
\end{equation*}

The general fiber of $\rho$ is rational since $\pi$ is locally
trivial with rational fibers. It follows that $\rho$ gives the MRC
fibration of $\M(a,e)$. Note that $\rho$ is dominant if and only
if $u$ is dominant, that is, when $\de\geq g$. \ep

\subsection{Geometric description of the varieties $\M(a,e)$
and $\M_E$}

\

\ps

\subsubsection{The varieties $\M(a,e)$}

Recall that we defined in \cite{C}, for every $\Lb\in\Pic^{-e}(C)$, a vector
space $V_{\Lb}$ as the space of extensions $\Ext^1_C(\Lb^{-1}\otimes\xi,\Lb)$.
Let $Z_{\Lb}\subset\PP(V_{\Lb})$ be the locus of unstable extensions.
We defined a morphism $$\ka_{\Lb}:\PP(V_{\Lb})\setm Z_{\Lb}\ra M$$
such that $\ka_{\Lb}^*\Th\cong O(2e+1)$. We recall that
$\dim\PP(V_{\Lb})=(2e+g-1)$.

\begin{lemma}\label{geom_descr}
If $f\in\M(a,e)$, then there is $\Lb\in\Pic^{-e}(C)$ such that $f$
is the extension of a rational map obtained as a composition
\begin{equation*}
\PP^1\dra\PP(V_{\Lb})\setm Z_{\Lb}\ra M
\end{equation*}

If $\de=0$, then the map  $\PP^1\ra\PP(V_{\Lb})\setm Z_{\Lb}$ is defined
everywhere and it has degree $(2a-k)$.
\end{lemma}

\bp
Let $f\in\M(a,e)$. The canonical  sequence of the vector bundle $\HH_f$ has
the form
\begin{equation}
0\ra O(a)\bo\Lb\ra\HH_f\ra(O(k-a)\bo\Lb^{-1}\otimes\xi)\otimes\I_Z\ra0.
\end{equation}

Let $\Ga\subset\PP^1$ be the set-theoretic image of $Z$ via the projection to
$\PP^1$. Recall that $Z$ is a $0$-cycle in $\PP^1\times C$ of length 
$\de=(k-a)-e(2a-k)$. The restriction of the morphism $f$ to 
$\PP^1\setm\Ga$ factors as
\begin{equation*}
\begin{CD}
\PP^1\setm\Ga@>g>>\PP(V_{\Lb})\setm Z_{\Lb}@>{\ka_{\Lb}}>> M
\end{CD}
\end{equation*}

If $\de=0$ then $Z=\eset$ and the morphism $g$ is
defined everywhere. As $f$ has degree $k$ and $\ka_{\Lb}^*\Th\cong
O(2e+1)$, it follows that $g$ has degree $\frac{k}{2e+1}=(2a-k)$.
\ep

Note that if $\de=0$, the variety $\Mae$ is the variety
$\S(e,n)$, for $n=(2a-k)$, defined in 
\cite{C} (Section 3.3). 
In particular, if $k$ is odd, then
$\overline{\M(\frac{k+1}{2},\frac{k-1}{2})}$
is the nice component of the space
$\Hom_k(\PP^1, M)$, as defined in \cite{C}.

\begin{lemma}\label{covering_M}
The only irreducible components $\M$ of the space $\Hom_k(\PP^1, M)$, 
of the form $\Mae$, for which the evaluation map $\ev:\PP^1\times\M\ra M$ 
is dominant, are for the following cases of $a$ and $e$ 
(see Figure 5): 
\bi
\item[i) ]If $k\geq(g-1)$ is an odd integer and
$a=\frac{k+1}{2}$, $e=\frac{k-1}{2}$
\item[ii) ]If $g$ is even and  $k$ is divisible by $(g-1)$ and
$a=\frac{kg}{2(g-1)}$, $e=\frac{g}{2}-1$
\ei
\end{lemma}
Note that by Lemma \ref{geom_descr}, in case i) the variety $\Mae$ is 
the nice component of \cite{C}, while in case ii) it is the 
\emph{almost nice component} of \cite{C}. 
The two components are the same if $g$ even and $k=(g-1)$.

\bp
By Lemma \ref{geom_descr}, we have that if $e<\frac{g}{2}-1$ 
there is a proper subvariety of $M$ that contains all the rational curves
in $\M$ (the images in $M$ of the morphisms $\ka_{\Lb}$ for all
$\Lb\in\Pic^{-e}(C)$ form a variety of dimension at most
$(2e+g-1)+g<(3g-3)$. The cases i) and ii) are precisely the irreducible 
components for which $e\geq\frac{g}{2}-1$. (See Figure 4 and the Figures in 
Section \ref{region_R}). 

To show that in these two cases, the evaluation map $\ev:\PP^1\times\M\ra M$ 
is dominant, we recall that we proved in 
\cite{C} (Corollary 4.2) 
that a general element of these components is a free curve 
$f:\PP^1\ra M$. It follows from the general theory of free rational curves
(see \cite{K}, II, 3.5.3) that these curves  cover $M$.
\ep

\subsubsection{The variety $\M_E$}

Assume that $k=2a$ is an even positive integer. In \cite{C}, for
every $D\in\Sym^a(C)$ and $\E\in M_{1-a}$ such that
$\det(\E)\cong\xi(-D)$, we defined a vector space $V_{D,\E}$ as the space 
of extensions $\Ext^1_C(O_D,\E)$. Let $Z_{D,\E}\subset\PP(V_{D,\E})$ be 
the locus of unstable extensions. We defined a morphism
\begin{equation*}
\eta_{D,\E}:\PP(V_{D,\E})\setm Z_{D,\E}\ra M
\end{equation*} such that
such that $\eta_{D,\E}^*\Th\cong O(k)$.
Recall that $\dim\PP(V_{D,\E})=k-1$.

\begin{lemma}
If  $f\in\M_E$, then there is  a pair $(D,\E)$ as before
such that $f$ is a composition
\begin{equation*}
\begin{CD}
\PP^1@>g>>\PP(V_{D,\E})\setm Z_{D,\E}@>{\eta_{D,\E}}>> M
\end{CD}
\end{equation*}
\end{lemma}

If $k$ is even, the variety $\ME$ is the nice component of the space
$\Hom_k(\PP^1, M)$, as defined in \cite{C}. By a similar argument as in
Lemma \ref{covering_M}, one may see that the rational curves in $\ME$
cover $M$, for any $k\geq2$. It is a consequence of the fact that
the general element of $\ME$ is a free curve for any even $k\geq2$
(see \cite{C}, Lemma 7.1). 

\subsection{Obstructed components}

\

\ps

Recall that an element $f\in\Hom_k(\PP^1,M)$ is unobstructed
if $h^1(\PP^1,f^*T_M)=0$. An unobstructed point is a smooth point
of $\Hom_k(\PP^1,M)$ and it is contained in a unique irreducible
component of the expected dimension $(2k+3g-3)$. Hence,
if an irreducible variety $\M\subset\Hom_k(\PP^1,M)$ has dimension
$\dim\M>(2k+3g-3)$, then $\M$ is obstructed.

\begin{prop}\label{obstr_comp}
Let $a$ and $e$ be integers such that $\M=\Mae$ is an
irreducible component of $\Hom_k(\PP^1,M)$. Then the general point
of $\M$ is unobstructed if and only if one of the following holds 
(see Figure 5): 
\bi \item[i) ]$k\geq(g-1)$ is an odd integer and $\M$ is the nice component:
$a=\frac{k+1}{2}$, $e=\frac{k-1}{2}$

\item[ii) ]$g$ is even, $k$ is divisible by $(g-1)$ and $\M$ is
the almost nice component: $a=\frac{kg}{2(g-1)}$,
$e=\frac{g}{2}-1$ \ei
\end{prop}

\begin{figure}[h]
\centerline{
\psfig{figure=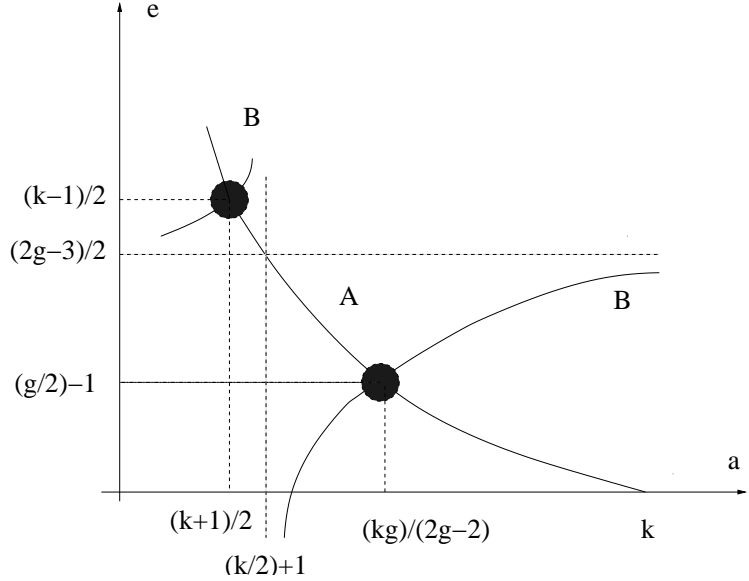,height=3in}
}
\smallskip
\centerline{
\parbox{4.5in}{\caption[]{{\small The nice component and the almost nice
component are the only irreducible components with unobstructed general
member }}} }
\end{figure}

\bp
Let $f\in\M(a,e)$ and, for simplicity, denote $\F=\HH_f$. We have
from \cite{N2} that $f^*T_{M}\cong\R^1{p_1}_*(\ad(\F))$, where
$\ad(\F)$ is the bundle of traceless endomorphisms. We have
\begin{equation*}
\H^1(\PP^1,\R^1{p_1}_*(\ad(\F))\cong\H^2(\PP^1\times C,\ad(\F))
\cong\H^2(\PP^1\times C,\F\otimes\F^*)
\end{equation*}

We compute $\H^2(\F\otimes\F^*)$ using the canonical sequence of $\F$:
\begin{equation}\label{can}
0\ra\F'\ra\F\ra\M\otimes\I_Z\ra0
\end{equation}

We first show that  $\H^2(\F'\otimes\F^*)=0$. By Serre duality,
if  $K$ is the canonical bundle of $\PP^1\times C$,
this is equivalent to $\H^0({\F'}^*\otimes\F\otimes K)=0$.
By tensoring (\ref{can}) with ${\F'}^*\otimes K$ and taking global sections,
we get an exact sequence:
\begin{equation*}
0\ra\H^0(K)\ra\H^0(\F\otimes{\F'}^*\otimes K)\ra
\H^0({\F'}^*\otimes K\otimes\M\otimes\I_Z)
\end{equation*}

Since $\H^0(K)=0$, it is enough to prove that
$\H^0({\F'}^*\otimes K\otimes\M\otimes\I_Z)=0$.
More generally, we show that ${\F'}^*\otimes K\otimes\M$ does
not have global sections. Note that
\begin{equation*}
\F'\cong O(a)\bo\Lb,\quad\M\cong O(k-a)\bo(\Lb^{-1}\otimes\xi)
\end{equation*}
where $\Lb$ is a line bundle on $C$ of degree $-e$. Since $K\cong
O(-2)\bo K_C$, it follows that
$${\F'}^*\otimes\K\otimes\M\cong O(k-2a-2)\bo(\Lb^{-2}\otimes\xi\otimes
K_C)$$ which does not have global sections, as $a>\frac{k}{2}$.
Hence, we have proved
\begin{equation*}
\H^2(\F'\otimes\F^*)\cong\H^0({\F'}^*\otimes\F\otimes K)=0
\end{equation*}

The long exact sequence coming from the sequence (\ref{can}) tensored with
$\F^*$ gives
\begin{equation*}
\H^2(\F\otimes\F^*)\cong\H^2(\M\otimes\I_Z\otimes\F^*)
\end{equation*}

By Serre duality, we have
$\H^2(\M\otimes\I_Z\otimes\F^*)\cong\H^0(K\otimes\M^*\otimes\F)$. By
tensoring the exact sequence (\ref{can}) with $K\otimes\M^*$, it
follows that
\begin{equation*}
\H^0(K\otimes\M^*\otimes\F)\cong\H^0(K\otimes\M^*\otimes\F')
\end{equation*}

The line bundle
\begin{equation*}
K\otimes\M^*\otimes\F\cong O(2a-k-2)\bo(\Lb^2\otimes\xi^{-1}\otimes K_C)
\end{equation*}
has global sections if $2a-k\geq2$ and $e<\frac{g}{2}-1$.

Note from Theorem \ref{main}, that the only pair of integers $a$ and $e$
outside this range for which $\Mae$ is an irreducible component are
the two mentioned cases. (See also Figure 5 and the Figures in Section
\ref{region_R}.) 
Since if $2a-k=1$ or if $e=\frac{g}{2}-1$ and $\Lb$ is general 
in $\Pic^{-e}(C)$ (use Brill-Noether theory, see 
\cite{ACGH}, IV.4.5) the line bundle $K\otimes\M^*\otimes\F$ has no global 
sections, it follows that the general point any of these irreducible
components is unobstructed.
\ep


\begin{cor}\label{non-red}
Let $\M$ be an irreducible component of the form $\M(a,e)$ that has
the expected dimension $(2k+3g-3)$ and it is different than the nice
component or the almost nice component. Then $\Hom_k(\PP^1,M)$ is
non-reduced along $\M$.
\end{cor}

\bp By Proposition \ref{obstr_comp}, a general point $f\in\M$ is
obstructed. We have $h^1(\PP^1, f^*T_M)\neq0$. It follows that
$h^0(\PP^1, f^*T_M)>(2k+3g-3)$. Since the dimension of $\M$ is
$(2k+3g-3)$ and the tangent space to $\Hom_k(\PP^1, M)$ at $f$ is
isomorphic to $H^0(\PP^1,f^*T_M)$, it follows that $\Hom_k(\PP^1, M)$ is
non-reduced at the point $f$. \ep

The pairs of integers $(a,e)$ in Corollary \ref{non-red} are
on the part of the hyperbola (\ref{B}) that
lies beneath the line $e=\frac{g}{2}-1$.
For example, if $g=10$ and $k=33$, the irreducible component
$\M(18,1)$ has the expected dimension, but the general point is obstructed.


\section{Particular cases and examples}

\textbf{The irreducible component with the largest dimension}

\ps

\begin{figure}[h]
\centerline{
\psfig{figure=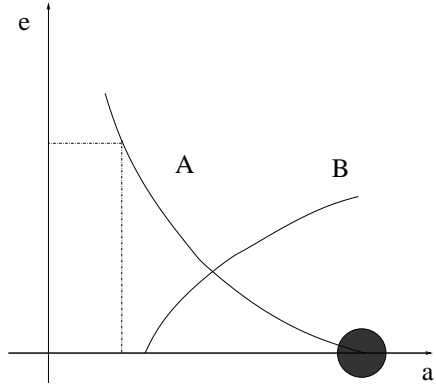,height=2in}}
\smallskip
\leftline{
\parbox{4.5in}{\caption[]{{\small
The component $M(k,0)$ has the largest dimension}}}
}
\end{figure}

Consider the variety $\M(k,0)$ (see Figure 6). 
 Using Proposition \ref{families_U},
one may prove that $\M(k,0)$ is closed in $\Hom_k(\PP^1,M)$.
By Theorem \ref{main}, the variety $\M(k,0)$ is an irreducible component of
$\Hom_k(\PP^1,M)$. It has dimension
\begin{equation*}
\dim\M(k,0)=(k+2)g-1
\end{equation*}

Using Observation \ref{dim_grow}, one may see that the variety 
$\M(k,0)$ is the irreducible component of the largest dimension.

By Lemma \ref{geom_descr}, a general element
$f\in\M(k,0)$ may be obtained from degree $k$ rational curves
in some $(g-1)$ dimensional linear subspaces of $M$. Note that
there is a $g$-dimensional family of such linear  subspaces of $M$, hence,
the rational curves in $\M(k,0)$ stay in a proper subvariety of $M$ of
dimension $(2g-1)$.


\

\textbf{Low degree rational curves}

\begin{cor}\textbf{(Lines)}
The variety $\Hom_1(\PP^1, M)$  has the
expected dimension $(3g-1)$ and is isomorphic
to a projective bundle over the variety $Pic^0(C)$.
\end{cor}

\begin{figure}[h]
\centerline{
\psfig{figure=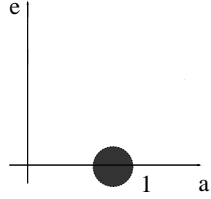,height=1in}
}
\smallskip
\centerline{
\parbox{4.5in}{\caption[]{{\small
If $k=1$, there is a unique irreducible component and it has the expected
dimension $(3g-1)$}}}
}
\end{figure}

\bp
Let $k=1$. The range (\ref{star}) for $a$ and $e$ is just the point
$(1,0)$ (see Figure 7). 
Let $\M$ be the  nice component $\M(1,0)$. It is the
unique  irreducible component of $\Hom_1(\PP^1, M)$ and it has
the expected dimension $(3g-1)$. By Corollary \ref{M_isom_B^0}, we have
$\M(1,0)\cong\Bb^0(1,0)$. It is easy to show that
$\Bb^0(1,0)=\Bb(1,0)=\Bb'(1,0)$ (see Theorem \ref{good_B(a,e)}).
It follows that $\M\cong\Bb'(1,0)$ is a projective bundle
over $\Pic^0(C)$.
\ep


\begin{cor}\textbf{(Conics)}
The scheme $\Hom_2(\PP^1,M)$ has two irreducible components:
the nice component $\ME$, which has the expected dimension $(3g+1)$ and the
component $\M(2,0)$ of dimension $(4g-1)$. An element in
$\M(2,0)$ is obtained by taking conics in some $(g-1)$-planes contained
in $M$.
\end{cor}

\bp
There is only one variety of the form $\M(a,e)$, namely $\M(2,0)$.
It follows from Theorem \ref{main} that $\M(2,0)$ and
$\ME$ are the only irreducible components. The geometric description
follows from Lemma \ref{geom_descr}.
\ep

\begin{cor}\textbf{(Cubics)}
The scheme $\Hom_3(\PP^1,M)$ has two irreducible components (see Figure 8):
the nice component $\ME$, which has the expected dimension
$(3g+3)$ and the component $\M(3,0)$ of dimension $(5g-1)$.
A general element in $\M(3,0)$ is obtained by taking cubics in
some $(g-1)$-planes contained in $M$.
\end{cor}

\bp

\begin{figure}[h]
\centerline{
\psfig{figure=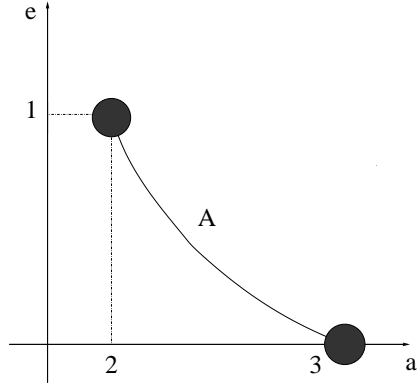,height=2in}
}
\smallskip
\centerline{
\parbox{4.5in}{\caption[]{{\small
If $k=3$, there are two irreducible components in the
space $\Hom_3(\PP^1,M)$}}}
}
\end{figure}

There are only two varieties of type $\M(a,e)$, namely $\M(2,1)$ and
$\M(3,0)$. It follows from Theorem \ref{main} that their closures
in $\Hom_3(\PP^1,M)$ are irreducible components. The geometric description
follows from Lemma \ref{geom_descr}.
\ep

\

\textbf{Low genus examples}

\

\begin{figure}[h]
\centerline{ \psfig{figure=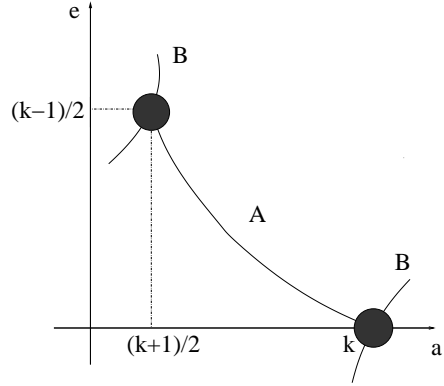,height=2in} }
\smallskip
\centerline{
\parbox{4.5in}{\caption[]{{\small
If $g=2$, the region $R$ consists of two points:
$\{(\frac{k+1}{2},\frac{k-1}{2}), (k,0)\}$}}}
}
\end{figure}

The following are consequences of Theorem \ref{main} and Lemma 
\ref{geom_descr}.
\begin{cor}\textbf{(Genus $2$)}
If $g=2$, for any $k\geq2$, there are only two irreducible components
in the space $\Hom_k(\PP^1, M)$, both of the expected dimension $(2k+3)$
(see Figure 9):
the nice component and the component $\M(k,0)$, whose elements are rational
curves $\PP^1\ra M$ which are $k$-to-$1$ onto lines contained in $M$.
\end{cor}
Note that the component $\M(k,0)$ is in this case the almost nice component.

\begin{figure}[h]
\centerline{ \psfig{figure=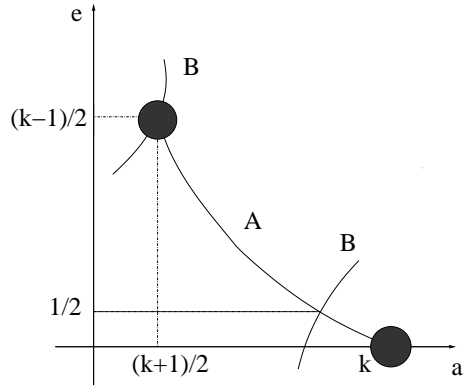,height=2in} }
\smallskip
\centerline{
\parbox{4.5in}{\caption[]{{\small
If $g=3$, the region $R$ consists of two points:
$\{(\frac{k+1}{2},\frac{k-1}{2}), (k,0)\}$}}}
}
\end{figure}

\begin{cor}\textbf{(Genus $3$)}
If $g=3$, for any $k\geq2$, there are only two irreducible components
in the space $\Hom_k(\PP^1, M)$ (see Figure 10) : 
the nice component, which has the expected dimension $(2k+6)$, and
the component $\M(k,0)$, of dimension $(3k+5)$, whose elements are
rational curves $\PP^1\ra M$ which are obtained by taking degree $k$ curves
in $2$-planes contained in $M$.
\end{cor}
Note that there is no almost nice component in this case.

\begin{figure}[h]
\centerline{ \psfig{figure=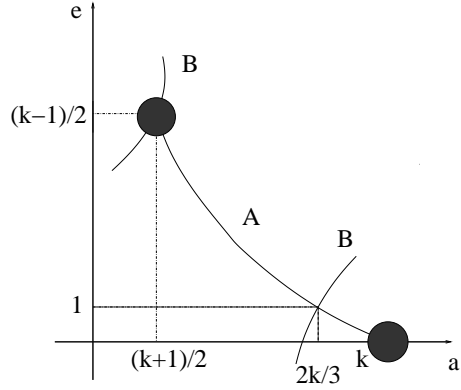,height=2in} }
\smallskip
\centerline{
\parbox{4.5in}{\caption[]{{\small
If $g=4$, the region $R$ consists of three points:
$\{(\frac{k+1}{2},\frac{k-1}{2}), (\frac{2k}{3},1), (k,0)\}$}}}
}
\end{figure}

\begin{cor}\textbf{(Genus $4$)}
If $g=4$, for all integers $k\geq2$ not divisible by $3$,
there are only two irreducible components in the space $\Hom_k(\PP^1, M)$
(see Figure 11): 
the nice component, which has the expected dimension $(2k+9)$, and
the component $\M(k,0)$, of dimension $(4k+7)$. If $k$ is divisible by $3$,
there is an extra component, namely, the almost nice component
$\overline{\M(\frac{2k}{3},1)}$, which also has the expected dimension.
\end{cor}

Note that in all the previous small genus cases, all the irreducible
components have $J(C)$ as MRC quotient. We give some examples where this
is not always the case.

\

\textbf{Example 1: $g=6$ and $k=15$ }

\

\begin{figure}[h]
\centerline{ \psfig{figure=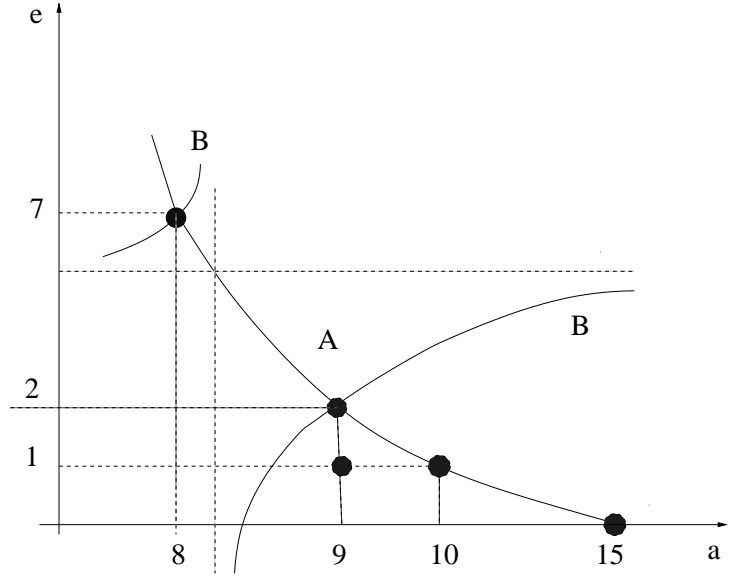,height=3in} }
\smallskip
\centerline{
\parbox{4.5in}{\caption[]{{\small There are five irreducible components
when $g=6$ and $k=15$}}}
}
\end{figure}

We have five irreducible components (see Figure 12): 
\bi
\item[i. ] The nice component $\overline{\M(8,7)}$, which has the expected
dimension  $45$
\item[ii. ] The almost nice component $\overline{\M(9,2)}$, which also has
dimension  $45$
\item[iii. ]The component $\overline{\M(9,1)}$, which has dimension $46$
\item[iv. ]The component $\overline{\M(10,1)}$, which has dimension $53$
\item[v. ]The component $\M(15,0)$, which has dimension $101$
\ei

\ps

Note that the MRC quotient of the component $\overline{\M(9,1)}$ is
isomorphic to $J(C)\times A$, where $A$ is the subvariety of $J(C)$
given by the image of the map $\Sym^3(C)\ra\Pic^3(C)$.
The MRC quotient of any of the other components is isomorphic to $J(C)$.

\

\textbf{Example 2: $g=7$ and $k=61$ }

\ps

\begin{figure}[h]
\centerline{ \psfig{figure=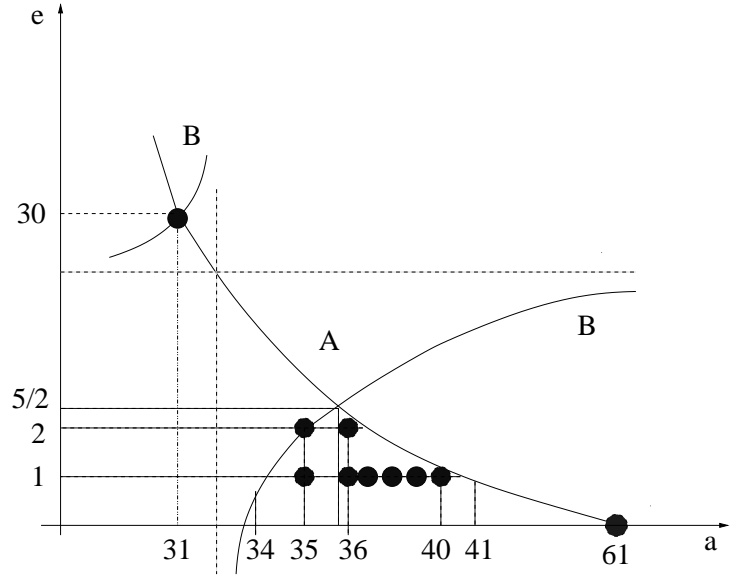,height=3in} }
\smallskip
\centerline{
\parbox{4.5in}{\caption[]{{\small There are $10$ irreducible components
when $g=7$ and $k=61$}}}
}
\end{figure}

In this case, we have $10$ irreducible components (see Figure 13):
\bi
\item[i. ] The nice component $\overline{\M(31,30)}$
\item[ii. ] The component $\M(61,0)$
\item[iii. ]The components $\overline{\M(35,2)}$ and $\overline{\M(36,2)}$
\item[iv. ]The components $\overline{\M(a,1)}$, for $a$ between
$35$ and $40$
\ei

\ps

Note that the MRC quotient of the components $\M(31,30)$ and
$\M(61,0)$ is $J(C)$, while the MRC quotient of all the other components
is not $J(C)$. If we compute $\de=(k-a)-e(2a-k)$, we see that
$\de\geq g$ in the following cases: $a=35, e=2$ and for
$35\leq a\leq 38, e=1$.
In all these cases the MRC quotient is $J(C)\times J(C)$.


\section{Pictures describing the region $R$}\label{region_R}

The following pictures describe the region $R$, corresponding to irreducible
components of the form $\Mae$.


\begin{figure}[h]
\centerline{ \psfig{figure=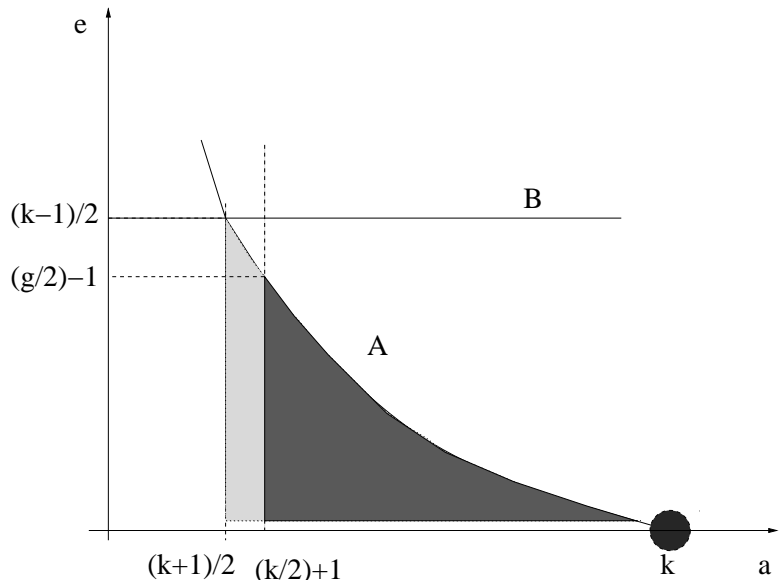,height=3in} }
\smallskip
\centerline{
\parbox{4.5in}{\caption[]{{\small
The dark shaded part is the region  $R$ when $g\geq3$ and $k=2g-2$}}}
}
\end{figure}

\begin{figure}[h]
\centerline{ \psfig{figure=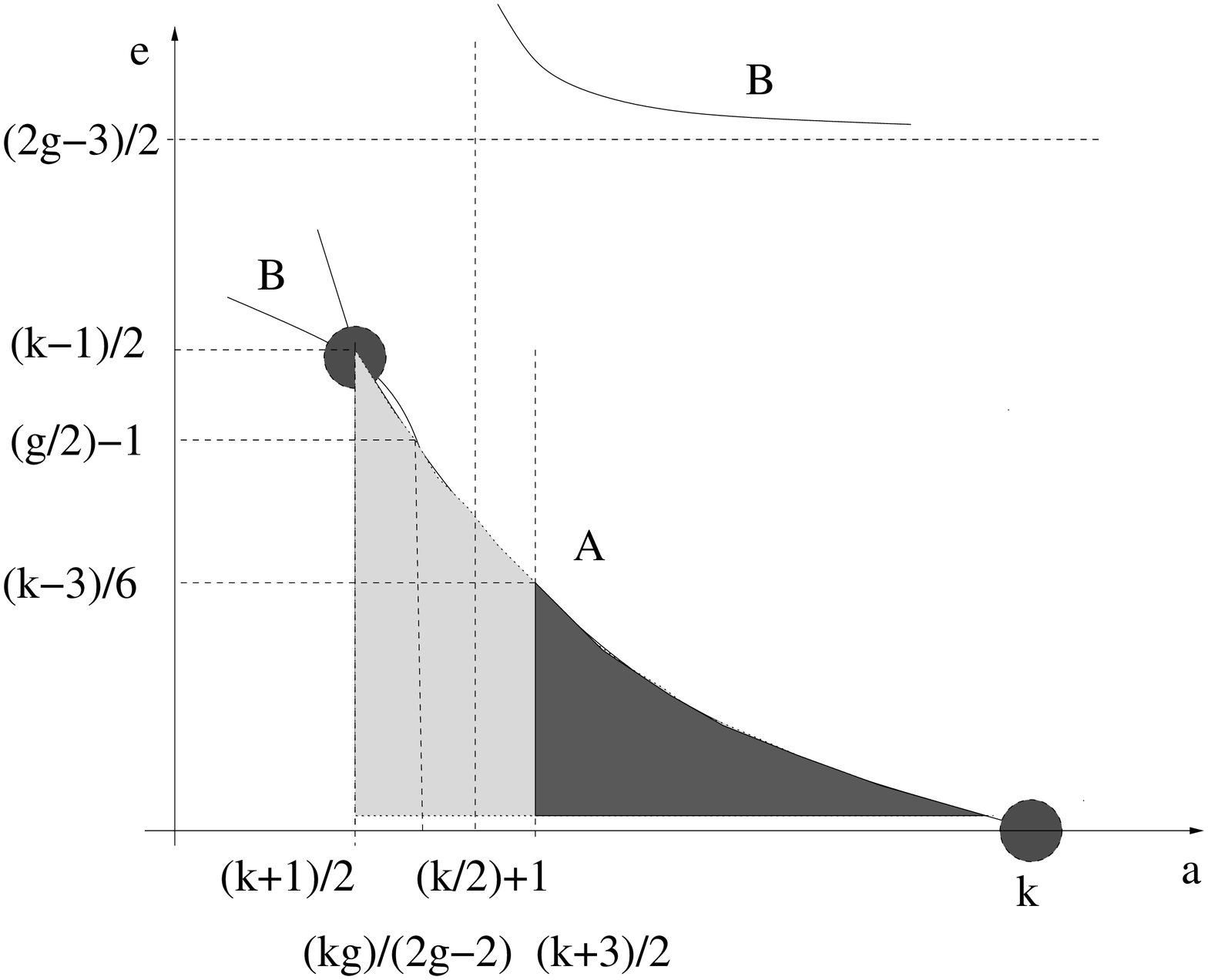,height=3in} }
\smallskip
\centerline{
\parbox{4.5in}{\caption[]{{\small
The dark shaded part is the region  $R$
when $g\geq3$ and $k$ is an odd integer in the interval
$g-1\leq k<2g-2$.}}}
}
\end{figure}

\begin{figure}[h]
\centerline{ \psfig{figure=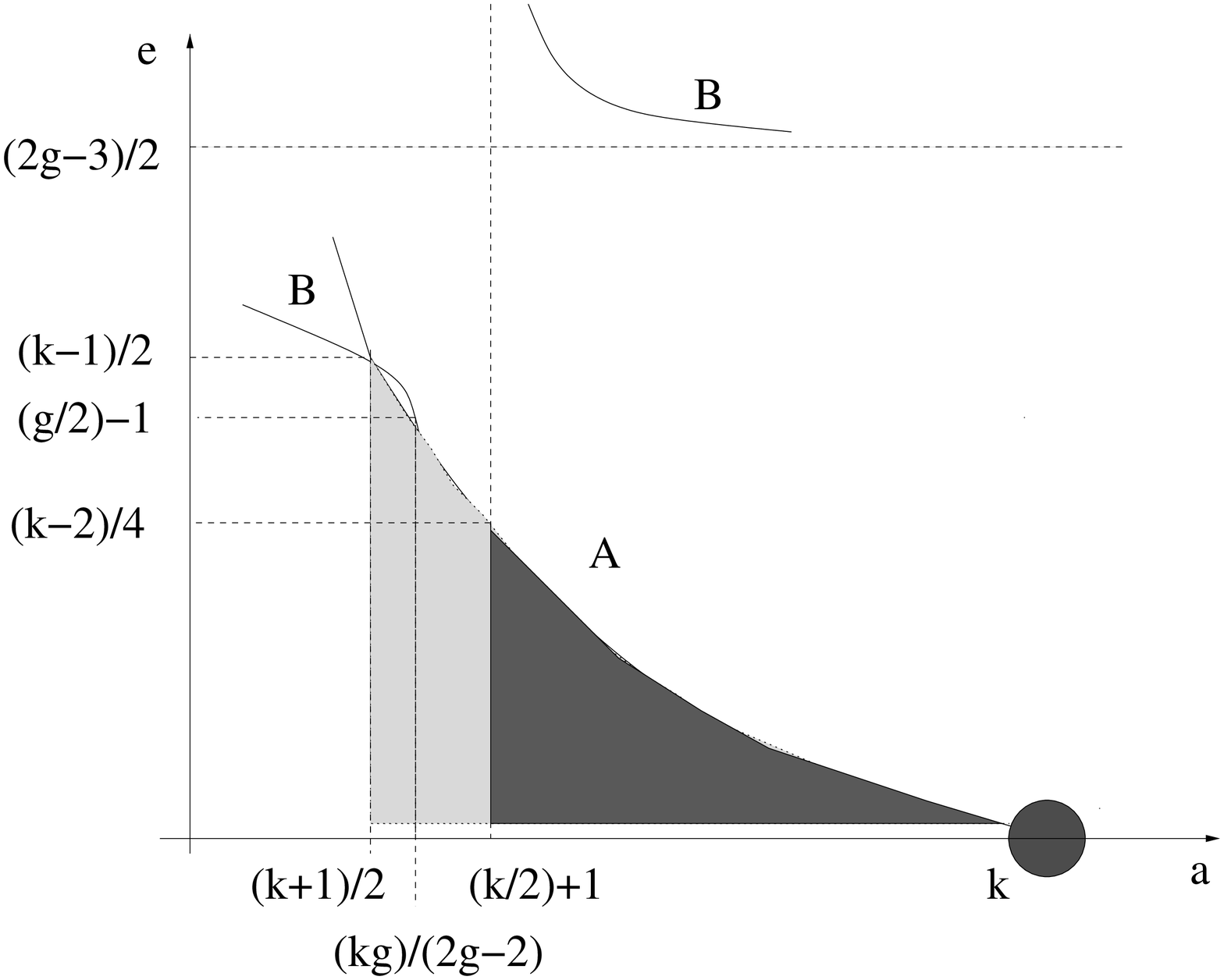,height=3in} }
\smallskip
\centerline{
\parbox{4.5in}{\caption[]{{\small
The dark shaded part is the region  $R$
when $g\geq3$ and $k$ is an even integer in the interval
$g-1\leq k<2g-2$}}}
}
\end{figure}

\

\begin{figure}[h]
\centerline{ \psfig{figure=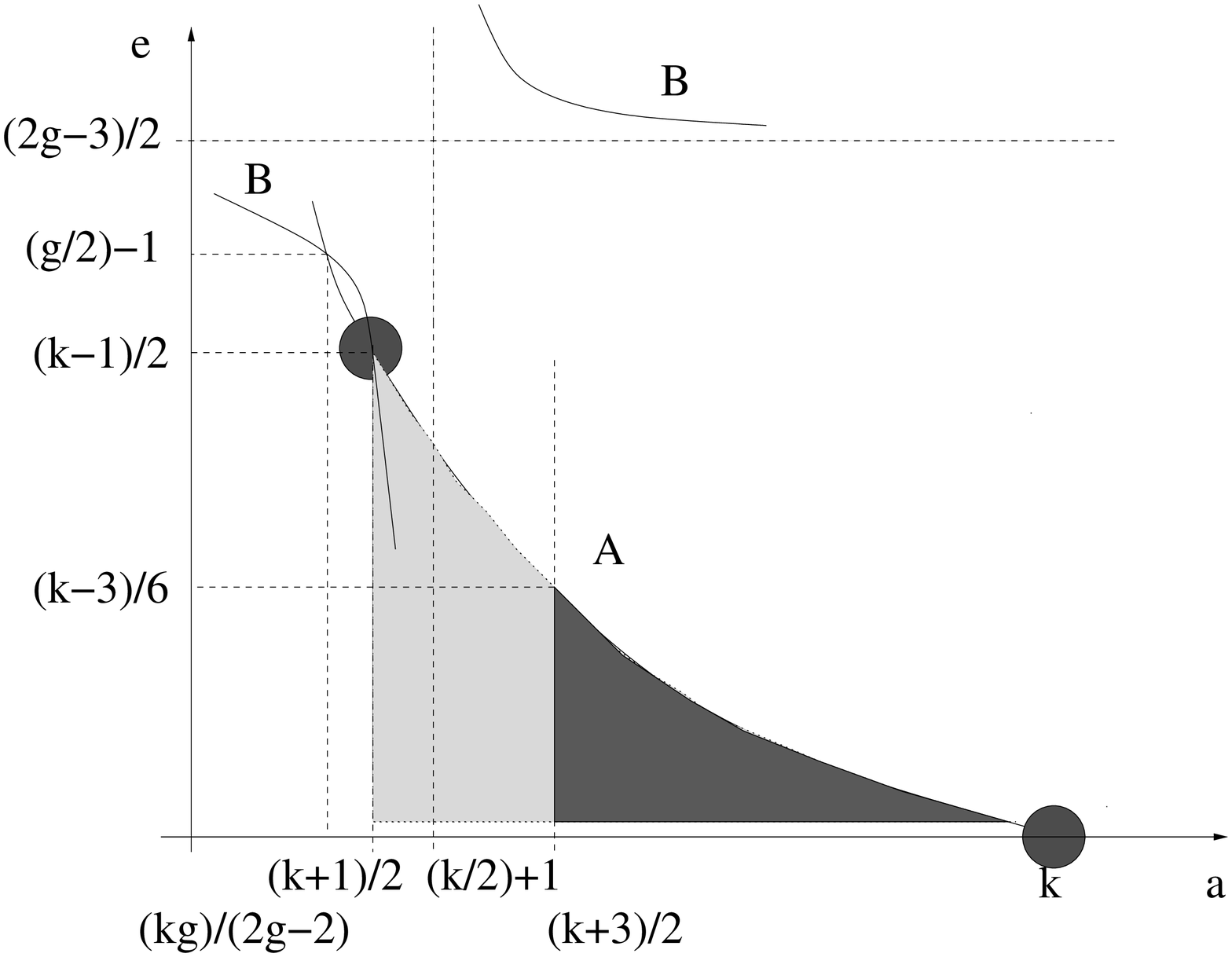,height=3in} }
\smallskip
\centerline{
\parbox{4.5in}{\caption[]{{\small
The dark shaded part is the region  $R$
when $g\geq3$ and $k$ is an odd integer in the interval $0<k<g-1$}}}
}
\end{figure}

\

\begin{figure}[h]
\centerline{ \psfig{figure=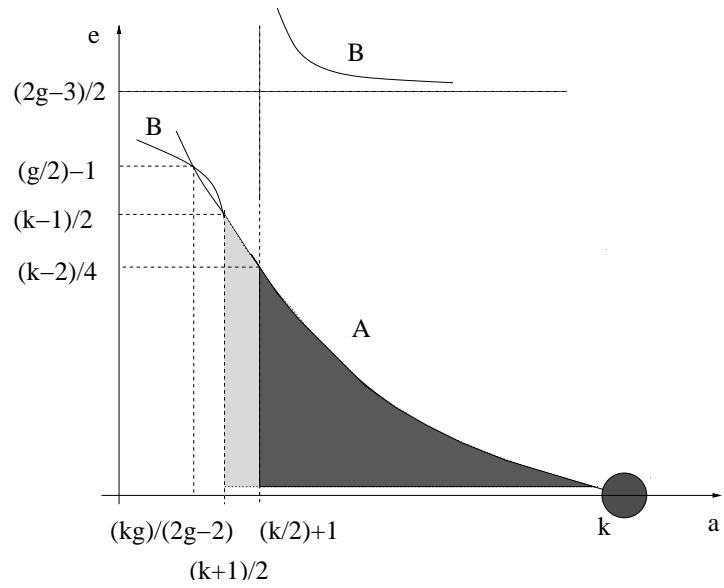,height=3in} }
\smallskip
\centerline{
\parbox{4.5in}{\caption[]{{\small
The dark shaded part is the region  $R$
when $g\geq3$ and $k$ is an even integer in the interval $0<k<g-1$}}}
}
\end{figure}




\begin{thebibliography}{xxxxxxxx}

\bibliographystyle{alpha}

\bibitem[ACGH]{ACGH}{E.~Arbarello, M. Cornalba, P.A. Griffiths, J. Harris},
{{\em Geometry of algebraic curves}},
{Springer-Verlag}{{}} {(1996)}


\bibitem[Brosius1]{BR1} {J.E. Brosius}, 
{{\em Rank $2$ Vector Bundles on a Ruled Surface. I}} ,
{Math. Ann.} 
{{\bf265}} {(1983)}, {} {155-168}


\bibitem[Brosius2]{BR2} {J.E. Brosius}, 
{{\em Rank $2$ Vector Bundles on a Ruled Surface. II}} ,
{Math. Ann.} 
{{\bf266}} {(1983)}, {} {199-214}


\bibitem[C1]{C} {A-M. Castravet},
{{Rational families of vector bundles I}} , {Pre-print}


\bibitem[DrezNar]{DN} {J.-M Drezet, M.S. Narasimhan}, 
{{\em Groupe de Picard des vari\'et\'es de modules
de fibr\'es semi-stables sur les courbes alg\'ebriques}},
{Invent. Math.} 
{{\bf97}} {(1989)}, {} {53-94}


\bibitem[Koll\'ar]{K} {J. Koll\'ar},
{{\em Rational Curves on Algebraic Varieties}},
{Springer-Verlag} {{}} {(1996)}



\bibitem[Mumford]{M} {D. Mumford},
{{\em Abelian Varieties}}, 
{Oxford University Press} {{}} {(1970)}




\bibitem[News72]{N2} {P.E. Newstead}, 
{{\em Characteristic classes of stable bundles 
of rank 2 over an algebraic curve}}, 
{Trans. Amer. Math. Soc. }
{{\bf169} } {(1972) }, {337-345 }



\bibitem[LePotier]{P} {J. Le Potier}, 
{{\em Lectures on Vector Bundles}}, 
{Cambridge Univ. Press} 
{{}} {(1997)} {} {}





\bibitem[Ramanan]{R} {S. Ramanan}, 
{{\em The Moduli Spaces of Vector Bundles over an 
Algebraic Curve}}, 
{Math. Ann.}
{{\bf200} } {(1973)}, {69--84}








\bibitem[Weibel]{W} {Charles A. Weibel}, 
{{\em An introduction to homological algebra}}, 
{Cambridge Univ. Press} 
{{}} {(1994)} {} {}



\end{thebibliography}
\end{document}